\documentclass[12pt,draftcls,onecolumn]{IEEEtran}

\ifCLASSINFOpdf
\else
\fi
\usepackage{graphicx}
\usepackage{epstopdf}
\usepackage{epsfig}
\usepackage{amsfonts}
\usepackage{amstext}
\usepackage{color}
\usepackage{empheq}
\usepackage{graphicx,colordvi}
\usepackage{amssymb,amsmath,amsopn,amsfonts,graphicx}
\usepackage{CJK}
\usepackage{float}

\newtheorem{lemma}{\bf Lemma}[section]
\newtheorem{theorem}{\bf Theorem}[section]
\newtheorem{definition}{\bf Definition}[section]

\newtheorem{remark}{\bf Remark}[section]
\newtheorem{corollary}{\bf Corollary}[section]
\newtheorem{proposition}{\bf Proposition}[section]
\newtheorem{example}{\bf Example}[section]

\usepackage[justification=centering]{caption}

\begin{document}
%
%
%
%

\title{Mean Field Linear Quadratic Control: FBSDE and Riccati Equation Approaches}

\author{
	Bingchang Wang, \emph{Member, IEEE,}
	\thanks{This work was supported by the National Natural Science Foundation of China under Grants 61773241, 61573221 and
61633014.}
	\thanks{Bingchang Wang is with the School of Control Science and Engineering, Shandong University, Jinan 250061, P. R. China. (e-mail: bcwang@sdu.edu.cn) }
and Huanshui Zhang, \emph{Senior Member, IEEE}

	\thanks{Huanshui Zhang is with the School of Control Science and Engineering, Shandong University, Jinan 250061, P. R. China. (e-mail: hszhang@sdu.edu.cn) }

}

%
%

\markboth{Journal of \LaTeX\ Class Files}%
{Shell \MakeLowercase{\textit{et al.}}: Bare Demo of IEEEtran.cls for IEEE Journals}
%



\maketitle

\begin{abstract}
 This paper studies social optima and Nash games for mean field linear quadratic control systems,
where subsystems are coupled via dynamics and individual costs. 
For the social control problem, we first obtain a set of forward-backward stochastic differential equations (FBSDE) from variational analysis, and construct a feedback-type control by decoupling the FBSDE. By using solutions of two Riccati equations, we design a set of decentralized control laws, which is further proved to be asymptotically social optimal. Two equivalent conditions are given for uniform stabilization of the systems in different cases.
For the game problem, we first design a set of decentralized control from variational analysis, and then show that such set of decentralized control constitute an asymptotic Nash equilibrium by exploiting the stabilizing solution of a nonsymmetric Riccati equation.

It is verified that the proposed decentralized control laws are equivalent to the feedback strategies of mean field control in previous works. This may illustrate the relationship between open-loop and feedback solutions of mean field control (games).
\end{abstract}


%
\IEEEpeerreviewmaketitle

\section{Introduction}
Mean field games 
have drawn increasing attention in many fields including system control, applied mathematics and economics \cite{BFY13, C14, GS13}. The mean field game involves a very large population of small interacting players with the feature that while the influence of each one is negligible, the impact of the overall population is significant. By combining mean field approximations and individual's best response,
the dimensionality difficulty is overcome. Mean field games and control
have found wide applications, including smart grids \cite{MCH13, CBM15}, finance, economics \cite{GLL11, CS14, WH15}, 
and social sciences 
\cite{BTN16}, etc.


By now, mean field games have been intensively studied in the LQ (linear-quadratic) framework \cite{HCM03, HCM07, LZ08, WZ13, BSYY16, MB17}.  
Huang \emph{et al.} developed the Nash certainty equivalence (NCE) based on the fixed-point method and designed an $\epsilon$-Nash equilibrium for mean field LQ games with discount costs by the NCE approach \cite{HCM03, HCM07}. The NCE approach was then applied to the cases with long run average costs \cite{LZ08} and with Markov jump parameters \cite{WZ13}, respectively. Bensoussan \emph{et al.} employed the adjoint equation approach and the fixed-point
theorem to obtain a sufficient condition for the unique existence of the
equilibrium strategy over a finite horizon \cite{BSYY16}.
For other aspects of mean field games, readers are referred to  \cite{HMC06, LL07, YMMS12, CD13} for nonlinear mean field games, \cite{weintraub2008markov} for oblivious equilibrium in dynamic games, \cite{H10, WZ12, WZ14} for mean field games with major players, \cite{HH16, MB17} for robust mean field games.

Besides noncooperative games, social optima in mean field models have also attracted much interest. The social optimum control refers to that all the players cooperate to optimize the common social cost---the sum of individual costs, which is usually regarded as a type of team decision problem \cite{R62, H80}. Huang \emph{et al.} considered social optima in mean field LQ control, and provided an asymptotic team-optimal solution \cite{HCM12}. Wang and Zhang \cite{WZ17} investigated a mean field social optimal problem where the Markov jump parameter appears as a common source of randomness. 
For further literature, see \cite{HN16} for social optima in mixed games, \cite{AM15} for team-optimal control with finite population and partial information. 

{Most previous results on mean field games and control were given by virtue of the fixed-point analysis. 
However, the fixed-point method is sometimes conservative, particularly for general systems. In this paper, we break away from the fixed-point method and solve the problem by tackling
forward-backward stochastic differential equations (FBSDE). In recent years, some substantial progress for the optimal LQ control has been made by solving the FBSDE. See \cite{Y13, ZX17, ZQ16, SLY16} for details.}

This paper investigates social optima and Nash games for linear quadratic mean field systems, 
where subsystems (agents) are coupled via dynamics and individual costs. 
For the finite-horizon social control problem, we first obtain a set of forward-backward stochastic differential equations (FBSDE) by examining the variation of the social cost, and give a centralized feedback-type control laws by decoupling the FBSDE. With mean field approximations, we design a set of decentralized control laws, which is further shown to have asymptotic social optimality.
For the infinite-horizon case, we design a set of decentralized control laws by using solutions of two Riccati equations, which is shown to be asymptotically social optimal. Some equivalent conditions are further given for uniform stabilization of the multiagent systems when the state weight $Q$ is semi-positive definite or only symmetric.
For the problem of mean field games, we first design a set of decentralized control by variational analysis, whose control gain satisfies a nonsymmetric Riccati equation. With the help of the stabilizing solution of the nonsymmetric Riccati equation, we show that the set of decentralized control laws is an asymptotic Nash equilibrium.
It is verified that the proposed decentralized control laws are equivalent representation of the feedback strategies in previous works of mean field control and games.
Finally, some numerical examples are given to illustrate the effectiveness of the proposed control laws.

The main contributions of the paper are summarized as follows.

(i) For the social control problem, we first obtain necessary and sufficient existence conditions of finite-horizon centralized optimal control by variational analysis, and then design a feedback-type decentralized control by tackling FBSDE with mean field approximations.

(ii) In the case $Q\geq0$, the necessary and sufficient conditions are given for uniform stabilization of the systems with the help of the system's observability and detectability.

(iii) In the case that $Q$ is only symmetric, the necessary and sufficient conditions are given for uniform stabilization of the systems using the Hamiltonian matrices.

(iv) For the game problem, we show that the decentralized control laws constitute an $\varepsilon$-Nash equilibrium by exploiting the stabilizing solution of a nonsymmetric Riccati equation.

(v) It is under nonconservative assumptions that we obtain the asymptotically optimal decentralized control, and such control laws are shown to be equivalent to the feedback strategies given by the fixed-point method in previous works \cite{HCM07, HCM12}.

The organization of the paper is as follows. In Section II, the socially optimal control problem is investigated. We first construct asymptotically optimal  decentralized control laws by tackling FBSDE for the finite-horizon case, then design asymptotically optimal control for the infinite-horizon case and further give two equivalent conditions of uniform stabilization for different cases. In Section III, we design a decentralized $\varepsilon$-Nash equilibrium for the finite-horizon and infinite-horizon cases, respectively. The proposed decentralized control laws are compared with the feedback strategies of previous works in Section IV. In Section V, some numerical examples are given to show the effectiveness of the proposed control laws. Section VI concludes the paper.

The following notation will be used throughout this paper. $\|\cdot\|$
denotes the Euclidean vector norm or matrix spectral norm. For a vector $z$ and a matrix $Q$, $\|z\|_Q^2= z^TQz$, and $Q>0$ ($Q\geq0$) means that $Q$ is positive definite (semi-positive definite). For two vectors $x,y$, $\langle x,y\rangle=x^Ty$.
$C([0,T],\mathbb{R}^n)$ is the space of all $\mathbb{R}^n$-valued continuous functions defined on $[0,T]$,
 and $C_{\rho/2}([0,\infty), \mathbb{R}^n)$ is a subspace of $C([0,\infty),\mathbb{R}^n)$ which is given by $\{f|\int_0^{\infty}e^{-\rho t}\|f(t)\|^2dt<\infty  \}.$
$L^2
_{\cal F}(0, T; \mathbb{R}^k)$ is the space of all $\mathcal{F}$-adapted $\mathbb{R}^k$-valued processes $x(\cdot)$ such that
$\mathbb{E}\int_0^T\|x(t)\|dt<\infty$.
For two sequences $\{a_n, n = 0, 1, \cdots\}$ and $\{b_n, n = 0, 1,  \cdots\}$, $a_n = O(b_n)$ denotes that
$\limsup_{n\to\infty}|{a_n}/{b_n}|\leq C$, and $a_n = o(b_n)$ denotes
$\limsup_{n\to\infty}|{a_n}/{b_n}|=0$.
For convenience of  presentation, we use $C, C_1,C_2,\cdots$ to
denote generic positive constants, which may vary from place to place.

\section{Mean Field LQ Social Control}\label{sec2.3.1}

Consider a large population systems with $N$ agents. Agent $i$ evolves by the following stochastic differential equation:
\begin{equation}\label{eq1}
\begin{aligned}
dx_i(t) = \ & [Ax_i(t)+Bu_i(t)+Gx^{(N)}(t)+f(t)]dt+\sigma(t) dW_i(t), \quad  1\leq i\leq N,
\end{aligned}
\end{equation}
where $x_i\in
\mathbb{R}^n$ and $u_i\in\mathbb{R}^r$ are the state and input of the $i$th agent.  $x^{(N)}(t)=\frac{1}{N}\sum_{j=1}^Nx_j(t)$, $f, \sigma\in C_{\rho/2}([0,\infty), \mathbb{R}^n)$.
$\{W_i(t),1\leq i\leq N\}$ are a sequence of independent $1$-dimensional Brownian motions on a complete
filtered probability space $(\Omega,
\mathcal F, \{\mathcal F_t\}_{0\leq t\leq T}, \mathbb{P})$.
The cost function of agent $i$ is given by
\begin{equation}\label{eq2}
\begin{aligned}
J_i(u)= \mathbb{E}\int_0^{\infty}e^{-\rho t}\Big\{&\big\|x_i(t)
-\Gamma x^{(N)}(t)-\eta\big\|^2_{Q}+\|u_i(t)\|^2_{R}\Big\}dt,
\end{aligned}
\end{equation}
where $Q$, $R$ are symmetric matrices with appropriate dimensions, and $R>0$.
Denote $u=\{u_1,
\ldots,u_{i}, \ldots, u_N\}$.
The decentralized control set is given by $$
\begin{aligned}
{\cal U}_{d,i} =\Big\{u_i\ \big|\ &
u_i(t)\ \hbox{is adapted to}\  \sigma(x_i(s),0\leq s\leq t),\ \mathbb{E}\int_0^{\infty}e^{-\rho t}\|u_i(t)\|^2dt<\infty\Big\}.
\end{aligned}
$$ For comparison, define the centralized control sets as
$$\begin{aligned}
  {\cal U}_{c,i} =\Big\{u_i \big|\ &
u_i(t)\ \hbox{is adapted to}\ {\mathcal F}_t, \ \mathbb{E}\int_0^{\infty}e^{-\rho t}\|u_i(t)\|^2dt<\infty\Big\},
\end{aligned}$$
and $
  {\cal U}_{c} =\big\{(u_1,\cdots,u_N) \big|\
u_i\ \hbox{is adapted to}\ {\cal U}_{c,i}$\big\},
where ${\mathcal F}_t= \sigma\{\bigcup_{i=1}^N{\mathcal F}_t^{i}\}$ and ${\mathcal F}_t^{i} = \sigma(x_{i}(0),\\ W_i(s),0\leq s\leq t), i=1,\cdots,N$.

In this section, we mainly study the following problem.

\textbf{(PS)}. Seek a set of decentralized control laws to optimize social cost
for the system (\ref{eq1})-(\ref{eq2}), i.e.,
$\inf_{u_i\in {\cal U}_{d,i}}J_{\rm soc},$
where $$J_{\rm soc}=\sum_{i=1}^NJ_i(u).$$

Assume

\textbf{A1)} $x_i(0), i=1,...,N$ are mutually independent and have the same mathematical expectation. $x_{i}(0)=x_{i0}$, $\mathbb{E}x_i(0)=\bar{x}_0$, $i=1,\cdots,N$. There exists a constant $C_0$ (independent of $N$) such that $\max_{1\leq i \leq N}\mathbb{E}\|x_i(0)\|^2<C_0$. Furthermore, $\{x_i(0), i=1,...,N\}$ and
$\{W_i, i=1,...,N\}$ are independent of each other.

\subsection{The finite-horizon problem}\label{sec3}
For the convenience of design, we first consider the following finite-horizon problem.
$$\textbf{(P1)} \inf_{u\in L^2_{{\cal F}_t}(0, T; \mathbb{R}^r)} J_{\rm soc}^{\rm F}(u),
$$
where $J_{\rm soc}^{\rm F}(u)=\sum_{i=1}^NJ_{i}^{\rm F}(u)$ and
\begin{equation}\label{eq3}
\begin{aligned}
J_{i}^{\rm F}(u)=\mathbb{E}\int_0^{T}e^{-\rho t}\Big\{&\big\|x_i(t)
-\Gamma x^{(N)}(t)-\eta\big\|^2_{Q}+\|u_i(t)\|^2_{R}\Big\}dt.
\end{aligned}
\end{equation}
We first give an equivalent condition for the convexity of Problem (P1).

\begin{proposition}\label{prop1}
 Problem (P1) is convex in $u$ if and only if
for any $u_i\in L^2_{{\cal F}_t}(0, T; \mathbb{R}^r)$, $i=1,\cdots,N$,
$$\sum_{i=1}^N\mathbb{E}\int_0^Te^{-\rho t}\Big\{\big\|y_i-\Gamma y^{(N)}\big\|^2_{Q}+\|u_i\|^2_{R}\Big\}dt\geq 0, $$
where $y^{(N)}=\sum_{j=1}^Ny_j/N$ and $y_i$ satisfies
\begin{align}\label{eq4aa}
  &dy_i=[A y_i+Gy^{(N)}+Bu_i]dt, \quad y_i(0)=0, \  i=1,2,\cdots,N.
  \end{align}
\end{proposition}

\emph{Proof.} Let $x_i$ and $\acute{x}_i$ be the state processes of agent $i$ 
with the control $v$ and $\acute{v}$,
respectively.
 Take any $\lambda_1\in [0, 1]$ and let $\lambda_2= 1-\lambda_1$.
Then
\begin{align*}
&\lambda_1  J_{\rm soc}^{\rm F} (v) +\lambda_2 J_{\rm soc}^{\rm F}  (\acute{v}) - J_{\rm soc}^{\rm F}(\lambda_1 v+\lambda_2 \acute{v})  \\
 =& \lambda_1\lambda_2 \mathbb{E}\int_0^T\left\{ \|x_i-\acute{x}_i -\Gamma(x^{N}-\acute{x}^{N})\|_Q^2 + \|u_i-\acute{u}_i\|^2_R\right\} dt.
\end{align*}
Denote $u= v-\acute{v}$, and $y_i=x_i-\acute{x}_i$. Thus, $y_i$ satisfies
(\ref{eq4aa}). 
By the definition of the convexity, the lemma follows.  $\hfill \Box$

By examining the variation of $ {J}_{\rm soc}^{\rm F}$, we obtain the necessary and sufficient conditions for
the existence of centralized optimal control of {(P1)}. 

\begin{theorem}\label{thm1}
	Suppose $R>0$. Then
(P1) has a
set of optimal control laws 
if and only if Problem (P1) is convex in $u$ and
	the following equation system admits a set of solutions $(x_i,p_i, \beta_{i}^{j},i,j=1,\cdots,N)$:
	\begin{equation}\label{eq4a}
	\left\{
	\begin{aligned}
	dx_i= &\big(Ax_i-B{R^{-1}}B^Tp_i+Gx^{(N)}+f\big)dt+\sigma dW_i,\\
	dp_i= &-\big[(A-\rho I)^Tp_i+G^Tp^{(N)}\big]dt-\big(Qx_i-Q_{\Gamma}x^{(N)}-\bar{\eta} \big)dt
	+\sum_{j=1}^N\beta_i^jdW_j,\\
	x_i(0)&={x_{i0}},\quad p_i(T)=0,\quad i=1,\cdots,N,
	\end{aligned}\right.
	\end{equation}
	where $Q_{\Gamma}\stackrel{\Delta}{=}\Gamma^TQ+Q\Gamma-\Gamma^TQ\Gamma$, $\bar{\eta}\stackrel{\Delta}{=}Q\eta-\Gamma^T Q\eta$, $p^{(N)}=\frac{1}{N}\sum_{i=1}^Np_i$,
	and
	furthermore the optimal control is given by $\check{u}_i=-{R^{-1}}B^Tp_i$.
\end{theorem}
{\it Proof.}
Suppose that $\check{u}_i=-R^{-1}B^Tp_i,$ where $p_i, i=1,\cdots,N$ are a set of solutions to the equation system
\begin{equation}\label{eq6}
\begin{aligned}
&dp_i=\alpha_idt+\beta_i^idW_i+\sum_{j\not =i}\beta_i^jdW_j, \quad p_i(T)=0, \quad i=1,\cdots, N,
\end{aligned}
\end{equation}
where $\alpha_i$, $i=1,\cdots,N$ are to be determined. 
 Denote by $\check{x}_i$ the state of agent $i$ under the control $\check{u}_i$. For any $u_i\in L^2_{{\cal F}_t}(0, T; \mathbb{R}^r) $ and $\theta\in \mathbb{R}$, let $u_i^{\theta}=\check{u}_i+\theta u_i$. Denote by $ x_i^{\theta}$ the solution of the following perturbed state equationㄩ
$$ \begin{aligned}
  &dx_i^{\theta}=\big(Ax_i^{\theta}+B(\check{u}_i+\theta u_i)+f+\frac{G}{N}\sum_{i=1}^Nx^{\theta}_i\big)dt+\sigma dW_i,\cr
 &x_i^{\theta}(0)=x_{i0},\ i=1,2,\cdots,N.
 \end{aligned}$$ Let $y_i=(x_i^{\theta}-\check{x}_i)/\theta$. 
It can be verified that
 $y_i$ satisfies (\ref{eq4aa}).
Then by It\^{o}'s formula, for any $i=1,\cdots,N$,
\begin{align*}
0
=\ &\mathbb{E}[\langle e^{-\rho T} p_i(T),y_i(T)\rangle-\langle p_i(0),y_i(0)\rangle]\cr
=\ &\mathbb{E}\int_0^T \big[\langle \alpha_i,y_i\rangle+\langle p_i,(A-\rho I)y_i+Gy^{(N)}+Bu_i\rangle\big] dt,
\end{align*}
which implies
\begin{align}\label{eq7}
0= \sum_{i=1}^N \mathbb{E}\int_0^T & e^{-\rho t} \big[\langle \alpha_i,y_i\rangle+\langle p_i,(A-\rho I)y_i+Gy^{(N)}+Bu_i\rangle\big] dt\cr
= \sum_{i=1}^N \mathbb{E}\int_0^T & e^{-\rho t} \Big[\langle \alpha_i+(A-\rho I)^Tp_i,y_i\rangle+\langle p_i,Bu_i\rangle\big] dt \cr &+\mathbb{E}\int_0^Te^{-\rho t}\big\langle \sum_{i=1}^Np_i,\frac{G}{N} \sum_{i=1}^Ny_i\big\rangle\Big] dt\cr
= \sum_{i=1}^N \mathbb{E}\int_0^T & e^{-\rho t} \big[\langle \alpha_i+(A-\rho I)^Tp_i+G^Tp^{(N)},y_i\rangle+\langle B^Tp_i,u_i\rangle\big] dt.
\end{align}
From (\ref{eq3}), we have
\begin{equation}\label{eq5a}
\begin{aligned}
 &\check{J}_{\rm soc}^{\rm F}(\check{u}+\theta u)-\check{J}_{\rm soc}^{\rm F}(\check{u})=2\theta I_1+{\theta^2}I_2
\end{aligned}
\end{equation}
where $\check{u}=(\check{u}_1,\cdots,\check{u}_N)$, and
\begin{align*}
I_1\stackrel{\Delta}{=}&\sum_{i=1}^N\mathbb{E}\int_0^Te^{-\rho t}  \big[\big\langle Q\big(\check{x}_i-(\Gamma\check{x}^{(N)}+\eta)\big),y_i-\Gamma y^{(N)}\big\rangle +
\langle R \check{u}_i,u_i\rangle \big]dt,\cr
I_2\stackrel{\Delta}{=}&\sum_{i=1}^N\mathbb{E}\int_0^Te^{-\rho t} \big[\big\|y_i
   -\Gamma y^{(N)}\big\|^2_{Q}
+\|u_i\|^2_{R}\big]dt.
\end{align*}
Note that
\begin{align*}
  &\sum_{i=1}^N\mathbb{E}\int_0^Te^{-\rho t} \big\langle Q\big(\check{x}_i-(\Gamma\check{x}^{(N)}+\eta)\big),\Gamma y^{(N)}\big\rangle dt\cr
  =&
\mathbb{ E}\int_0^Te^{-\rho t} \Big\langle \Gamma^TQ  \sum_{i=1}^N\big(\check{x}_i-(\Gamma\check{x}^{(N)}+\eta)\big),\frac{1}{N}  \sum_{j=1}^Ny_j\Big\rangle  dt\cr
 =& \sum_{j=1}^N \mathbb{E}\int_0^Te^{-\rho t} \Big\langle  \frac{\Gamma^TQ}{N} \sum_{i=1}^N\big(\check{x}_i-(\Gamma\check{x}^{(N)}+\eta)\big), y_j\Big\rangle  dt\cr
 =& \sum_{j=1}^N \mathbb{E}\int_0^Te^{-\rho t} \big\langle  {\Gamma^TQ} \big((I-\Gamma)\check{x}^{(N)}-\eta\big), y_j\big\rangle  dt.
\end{align*}
From (\ref{eq7}), one can obtain that
\begin{align}\label{eq10b}
 I_1=&\mathbb{E}\sum_{i=1}^N\int_0^Te^{-\rho t}\Big[\big\langle Q\big(\check{x}_i-(\Gamma\check{x}^{(N)}+\eta)\big),y_i-\Gamma y^{(N)}\big\rangle dt +
\langle R\check{u}_i+ B^Tp_i,u_i\rangle \Big]dt\cr
&+\sum_{i=1}^N \mathbb{E}\int_0^Te^{-\rho t} \big[\langle \alpha_i+(A-\rho I)^Tp_i+G^Tp^{(N)},y_i\rangle\cr
=&\sum_{i=1}^N\mathbb{E}\int_0^Te^{-\rho t}\Big\langle R\check{u}_i+B^Tp_i,u_i\Big\rangle dt\cr
&+\sum_{i=1}^N\mathbb{E}\int_0^Te^{-\rho t}\Big\langle
Q\big(\check{x}_i-(\Gamma\check{x}^{(N)}+\eta)\big)-{\Gamma^TQ} \big((I-\Gamma)\check{x}^{(N)}-\eta\big)\cr
&\quad +\alpha_i+(A-\rho I)^Tp_i+G^Tp^{(N)}, y_i\Big\rangle dt.
\end{align}
From (\ref{eq5a}), $\check{u}$ is a minimizer to Problem (P1) if and only if
$I_2\geq0$ and $I_1=0 $.
By Proposition \ref{prop1}, $I_2\geq0$ if and only if (P1) is convex. $I_1=0 $ is equivalent to
\begin{align*}
\alpha_i=&-\big[(A-\rho I)^Tp_i-\Gamma^TQ \big((I-\Gamma)\check{x}^{(N)}-\eta\big)+Q\big(\check{x}_i-(\Gamma\check{x}^{(N)}+\eta))+G^Tp^{(N)}\big],\cr
\check{u}_i=&-{R^{-1}}B^Tp_i.
\end{align*}
Thus, we have the following optimality system:
\begin{equation}\label{eq8}
\left\{
\begin{aligned}
d\check{x}_i=\ &(A\check{x}_i-B{R^{-1}}B^T\check{p}_i+G\check{x}^{(N)}+f)dt+\sigma dW_i,\\
d\check{p}_i=\ &-[(A-\rho I)^T\check{p}_i+G^T\check{p}^{(N)}+Q\check{x}_i-Q_{\Gamma}\check{x}^{(N)}+\bar{\eta} )]t+\sum_{j=1}^N\beta_i^jdW_j,\\
\check{x}_i(0)&={x_{i0}},\quad \check{p}_i(T)=0,\quad i=1,\cdots,N,
\end{aligned}\right.
\end{equation}
such that $\check{u}_i=-{R^{-1}}B^T\check{p}_i$. This implies that the equation systems (\ref{eq4a}) admits a solution\\ $(\check{x}_i,\check{p}_i,\check{\beta}_{i}^{j}, i,j=1,\cdots,N)$.

On other hand, if the equation system (\ref{eq4a}) admits a solution $(\check{x}_i,\check{p}_i, \check{\beta}_{i}^{j}, i,j =1,\cdots,N)$.
Let $\check{u}_i=-R^{-1}B^T\check{p}_i$. If (P1) is convex, then $\check{u}$ is a minimizer to Problem (P1).


$\hfill \Box$


It follows from (\ref{eq4a}) that
\begin{equation}\label{eq10}
\left\{
\begin{aligned}
d{x}^{(N)}=&\Big[(A+G){x}^{(N)}-B{R^{-1}}B^T{p}^{(N)}+f\Big]dt+\frac{1}{N}\sum_{i=1}^N\sigma dW_i,\\
d{p}^{(N)}=&-\Big[(A+G-\rho I)^T{p}^{(N)}-(I-\Gamma)^TQ(I-\Gamma){x}^{(N)}+\bar{\eta}\Big]dt+\frac{1}{N}\sum_{i=1}^N\sum_{j=1}^N{\beta}_i^jdW_j,\\
{x}^{(N)}(0)&=\frac{1}{N}\sum_{i=1}^Nx_{i0},\quad {p}^{(N)}(T)=0.
\end{aligned}\right.
\end{equation}
Let $p_i=Px_i+Kx^{(N)}+s$. Then by (\ref{eq4a}), (\ref{eq10}) and It\^{o}'s formula,
$$\begin{aligned}
dp_i= &P\Big[\big(A{x}_i-B{R^{-1}}B^T(Px_i+Kx^{(N)}+s)+Gx^{(N)}+f\big)dt+\sigma dW_i\Big]+(\dot{P}{x}_i+\dot{s}+\dot{K}x^{(N)})dt\\
&+K\Big\{\big[(A+G){x}^{(N)}-B{R^{-1}}B^T((P+K)x^{(N)}+s)+f\big]dt+\frac{1}{N}\sum_{i=1}^N\sigma dW_i\Big\}\cr
=\ &-\big[(A-\rho I)^T(Px_i+Kx^{(N)}+s)+G^T((P+K)x^{(N)}+s)+Qx_i-Q_{\Gamma}x^{(N)}-\bar{\eta} \big]dt\\
&+\sum_{j=1}^N\beta_i^jdW_j.
\end{aligned}$$
This implies that $\beta_i^i=\frac{1}{N}K\sigma+P\sigma$, $\beta_i^j=\frac{1}{N}K\sigma, \ j\not=i$,
\begin{align}\label{eq8a}
\rho{P}=\ &\dot{P}+A^TP+PA-PBR^{-1}B^TP+Q, \ P(T)=0,\\
\label{eq9a}
\rho K=\ &\dot{K}+(A+G)^TK+K(A+G)-PBR^{-1}B^TK-KBR^{-1}B^TP\cr
&+G^TP+PG-KBR^{-1}B^TK-Q_{\Gamma},\ K(T)=0,
\\ \label{eq10a}
\rho s=\ &\dot{s}+[A+G-BR^{-1}B^T(P+K)]^Ts+(P+K)f-\bar{\eta}, \quad s(T)=0.
\end{align}
Then $\check{u}_i=-{R^{-1}}B^T(Px_i+Kx^{(N)}+s).$

\begin{theorem}\label{thm2}
	Assume that A1) holds and $Q\geq0$. Then Problem (P1) has an optimal control $$\check{u}_i=-{R^{-1}}B^T(Px_i+Kx^{(N)}+s),$$
	where $P, K $ and $s$ are determined by (\ref{eq8a})-(\ref{eq10a}).
\end{theorem}
\emph{Proof.} Denote $\Pi=P+K$. Then from (\ref{eq9a}) and (\ref{eq10a}), $\Pi$ satisfies
\begin{equation}\label{eq11}
\begin{aligned}
\rho\Pi=\dot{\Pi}&+(A+G)^T\Pi+\Pi (A+G)-\Pi BR^{-1}B^T\Pi+\hat{Q}, \quad \Pi(T)=0,
\end{aligned}
\end{equation}
where $\hat{Q}\stackrel{\Delta}{=}(I-\Gamma)^TQ(I-\Gamma)$.
Note that $Q\geq0$ and $R>0$. By \cite{AM90, YZ99}, (\ref{eq8a}) and (\ref{eq11}) admit unique solutions $P\geq 0$
and $\Pi\geq 0$, respectively, which implies that (\ref{eq9a}) and (\ref{eq10a}) have unique solutions $K$ and $s$, respectively.
Then by \cite{MY99, ZX17}, the FBSDE (\ref{eq4a}) admits a unique solution. 
By Theorem \ref{thm1}, Problem (P1) has an optimal control
given by
$\check{u}_i=-{R^{-1}}B^T(Px_i+Kx^{(N)}+s),$
where $P, K $ and $s$ are determined by (\ref{eq8a})-(\ref{eq10a}).  \rightline{$\Box$}

As an approximation to ${x}^{(N)}$ in (\ref{eq10}), we obtain
\begin{equation}\label{eq12a}
\frac{d\bar{x}}{dt}=(A+G)\bar{x}-B{R^{-1}}B^T(\Pi\bar{x}+s)+f,\  \bar{x}(0)=\bar{x}_0.
\end{equation}
Then, by Theorem \ref{thm2}, the decentralized control law for agent $i$ may be taken as
\begin{equation}\label{eq15a}
\begin{aligned}
\hat{u}_i(t)=&-{R^{-1}}B^T(P\hat{x}_i(t)+K\bar{x}(t)+s(t)),\ 0\leq t\leq T,\ i=1,\cdots, N,
\end{aligned}
\end{equation}
where $P, K$, and $s$ are determined by (\ref{eq8a})-(\ref{eq10a}), and $\bar{x}$ and $\hat{x}_i$ satisfy (\ref{eq12a}) and
%
%
%
\begin{equation}\label{eq20}
d\hat{x}_i=\big[(A-BR^{-1}B^TP)\hat{x}_i-BR^{-1}B^T[K\bar{x}+s]+G\hat{x}^{(N)}+f\big]dt+\sigma dW_i.
\end{equation}

\begin{remark}
{ In previous works \cite{HCM12, WZ17}, the mean field term $x^{(N)}$ in cost functions (dynamics) is first substituted by a deterministic function $\bar{x}$. By solving an optimal tracking problem subject to consistency requirements, a fixed-point equation is obtained. The decentralized control is constructed by handling the fixed-point equation. Here, we firstly obtain the centralized open-loop solution by variational analysis.
  By tackling the coupled FBSDEs combined with mean field approximations, the decentralized control laws are designed. Note that in this case $s$ and $\bar{x}$ are fully decoupled and no fixed-point equation is needed.}
\end{remark}


\begin{theorem}\label{thm3}
	Let A1) hold and $Q\geq0$. 
The set of decentralized control laws
	$\{\hat{u}_1,\cdots,\hat{u}_N\}$ given by (\ref{eq15a}) has asymptotic social optimality, i.e.,
	$$\Big|\frac{1}{N}J^{\rm F}_{\rm soc}(\hat{u})-\frac{1}{N}\inf_{u\in L^2_{{\cal F}_t}(0, T; \mathbb{R}^{nr})}J^{\rm F}_{\rm soc}(u)\Big|=O(\frac{1}{\sqrt{N}}).$$
	
\end{theorem}

\emph{Proof.} See Appendix A. $\hfill \Box$

\subsection{The infinite-horizon problem}

Based on the analysis in Section \ref{sec3}, we may design the following decentralized control laws for Problem (PS):
\begin{equation}\label{eq14}
\begin{aligned}
\hat{u}_i(t)=\ &-{R^{-1}}B^T[P\hat{x}_i(t)+(\Pi-P)\bar{x}(t)+s(t)],\ t\geq 0, \ \ i=1,\cdots, N,
\end{aligned}
\end{equation}
where $P$ and $\Pi$ are determined by
\begin{align}\label{eq15}
\rho P=&A^TP+PA-PBR^{-1}B^TP+Q,\\
\label{eq16}
\rho \Pi= &(A+G)^T\Pi+\Pi (A+G)-\Pi BR^{-1}B^T\Pi+\hat{Q},
\end{align}
and
$s, \bar{x}\in C_{\rho/2}([0,\infty),\mathbb{R}^n)$ are determined by
\begin{align}
\rho s&=\dot{s}+[A+G-BR^{-1}B^T\Pi]^Ts+\Pi f-\bar{\eta},\label{eq17}\\
\frac{d\bar{x}}{dt}&=(A+G)\bar{x}-B{R^{-1}}B^T(\Pi\bar{x}+s)+f,\   \bar{x}(0)=\bar{x}_0.\label{eq18}
\end{align}
Here the existence conditions of $P, \Pi,s$ and $\bar{x}$ need to be investigated further.

We introduce some assumptions:

\textbf{A2)} The system $(A-\frac{\rho}{2} I, B)$ is stabilizable, and $(A+G-\frac{\rho}{2} I, B)$ is stabilizable. 

\textbf{A3)}   $Q\geq0$, $(A-\frac{\rho}{2}I, \sqrt{Q}$) is observable, and 
$(A+G-\frac{\rho}{2}I, \sqrt{Q}(I-\Gamma))$ is observable. 

Assumptions A2) and A3) are basic in the study of the LQ optimal control problem. We will show that under some conditions, A2) is also necessary for uniform stabilization of multiagent systems.
In many cases, A3) may be weakened to the following assumption. 

{\bf{A3$^{\prime}$)}} $Q\geq0$, $(A-\frac{\rho}{2}I, \sqrt{Q}$) is detectable, and 
$(A+G-\frac{\rho}{2}I, \sqrt{Q}(I-\Gamma))$ is detectable. 


\begin{lemma}\label{lem2a}
	Under A2)-A3), (\ref{eq15}) and (\ref{eq16}) admit unique solutions $P>0, \Pi>0$, respectively, and (\ref{eq17})-(\ref{eq18}) admits a set of unique solutions $s, \bar{x}\in C_{\rho/2}([0,\infty),\mathbb{R}^n)$.
\end{lemma}

\emph{Proof.} From A2)-A3) and \cite{AM90}, (\ref{eq15}) and (\ref{eq16}) admit unique solutions $P>0, \Pi>0$ such that
$A-BR^{-1}B^TP-\frac{\rho}{2}I$ and $A+G-BR^{-1}B^T\Pi-\frac{\rho}{2}I$ are Hurwitz, respectively. From an argument in \cite[Appendix A]{WZ12}, we obtain
$s\in C_{\rho/2}([0,\infty),\mathbb{R}^n)$ if and only if $$s(0)=\int_0^{\infty}e^{(A+G-BR^{-1}B^T\Pi-{\rho}I)\tau}(\Pi f-\bar{\eta})d\tau.
$$
Under this initial condition, we have
$$s(t)=\int_t^{\infty}e^{-(A+G-BR^{-1}B^T\Pi-{\rho}I)(t-\tau)}(\Pi f-\bar{\eta})d\tau.$$
It is straightforward that $\bar{x}\in C_{\rho/2}([0,\infty),\mathbb{R}^n)$.
$\hfill \Box$


We further introduce the following assumption.

\textbf{A4)} $\bar{A}+G-\frac{\rho}{2}I$ is Hurwitz, where $\bar{A}\stackrel{\Delta}{=}A-BR^{-1}B^TP$.

\begin{lemma}\label{lem2}
	Let A1)-A4) hold. Then for (PS),
	\begin{equation}\label{eq13a}
	\mathbb{E}\int_0^{\infty} e^{-\rho t}\|\hat{x}^{(N)}(t)-\bar{x}(t)\|^2dt=O(\frac{1}{N}).
	\end{equation}

\end{lemma}
\emph{Proof.} See Appendix \ref{app b}.  $\hfill \Box$

It is shown that the decentralized control laws (\ref{eq15a}) uniformly stabilize the systems (\ref{eq1}) .

\begin{theorem}\label{thm4}
	Let A1)-A4) hold.  Then for any $N$, 
	\begin{equation}\label{eq13b}
	\sum_{i=1}^N\mathbb{E}\int_0^{\infty} e^{-\rho t} \left(\|\hat{x}_i(t)\|^2+\|\hat{u}_i(t)\|^2\right)dt<\infty.
	\end{equation}

\end{theorem}

\emph{Proof.} See Appendix \ref{app b}.  $\hfill \Box$

We now give two equivalent conditions for uniform stabilization of multiagent systems.

\begin{theorem}\label{thm5}
	Let A3) hold. Then for (PS) the following statements are equivalent:
	
	(i)  For any initial condition $(\hat{x}_1(0),\cdots, \hat{x}_N(0))$ satisfying A1),
	\begin{equation}\label{eq23}
	\sum_{i=1}^N\mathbb{E}\int_0^{\infty} e^{-\rho t} \left(\|\hat{x}_i(t)\|^2+\|\hat{u}_i(t)\|^2\right)dt<\infty.
	\end{equation}
	
	(ii) (\ref{eq15}) and (\ref{eq16}) admit unique solutions $P>0, \Pi>0$, respectively, and $\bar{A}+G-\frac{\rho}{2}I$ is Hurwitz.
	
	(iii) A2) and A4) hold.
	
\end{theorem}

\emph{Proof.} See the Appendix C.  $\hfill \Box$


For the case $G=0$, we have a simplified version of Theorem \ref{thm5}.
\begin{corollary}
	Assume that A3) holds and $G=0$. Then for (PS) the following statements are equivalent:
	
	(i)  For any $(\hat{x}_1(0),\cdots, \hat{x}_N(0))$ satisfying A1),
	\begin{equation*}
	\sum_{i=1}^N\mathbb{E}\int_0^{\infty} e^{-\rho t} \left(\|\hat{x}_i(t)\|^2+\|\hat{u}_i(t)\|^2\right)dt<\infty.
	\end{equation*}	
	
	(ii) (\ref{eq15}) and (\ref{eq16}) admit unique solutions $P>0, \Pi>0$, respectively.
	
	(iii) A2) holds.
\end{corollary}

When A3) is weakened to A3$^{\prime}$), we have the following equivalent conditions of uniform stabilization of the systems.
\begin{theorem}\label{thm6}
	Let A3$^{\prime}$) hold. Then for (PS) the following statements are equivalent:
	
	(i)  For any initial condition $(\hat{x}_1(0),\cdots, \hat{x}_N(0))$ satisfying A1),
	\begin{equation*}
	\sum_{i=1}^N\mathbb{E}\int_0^{\infty} e^{-\rho t} \left(\|\hat{x}_i(t)\|^2+\|\hat{u}_i(t)\|^2\right)dt<\infty.
	\end{equation*}
	
	(ii) (\ref{eq15}) and (\ref{eq16}) admit unique solutions $P\geq0, \Pi\geq0$, respectively, and $\bar{A}+G-\frac{\rho}{2}I$ is Hurwitz.
	
	(iii) A2) and A4) hold.
	
\end{theorem}

\emph{Proof.} See the Appendix C.  $\hfill \Box$

\begin{remark}
 { In \cite{ZQ16}, some similar results were given for the stabilization of mean field systems. However, only the limiting problem
  was considered in their work and the mean field term in dynamics and costs is $\mathbb{E}x(t)$ instead of $x^{(N)}$. Here we study large-population multiagent systems
  and the number of agents is large but not infinite. 
  The errors of mean field approximations are further analyzed.
   To obtain asymptotic optimality, an additional assumption
  A4) is needed later.}
\end{remark}

For the more general case that $Q$ are only symmetric, we have the following equivalent
conditions for uniform stabilization of multiagent systems.

Denote
$$M_1=\left[\begin{array}{cc}
A-\frac{\rho}{2}I & BR^{-1}B^T \\
Q& -A^{T}+\frac{\rho}{2}I
\end{array}\right],\quad
M_2=\left[\begin{array}{cc}
A+G-\frac{\rho}{2}I & BR^{-1}B^T \\
\hat{Q}& -(A+G)^{T}+\frac{\rho}{2}I
\end{array}\right].$$

\begin{theorem}\label{thm7}
	Assume that both $M_1$ and $M_2$ have no eigenvalues on the imaginary axis. Then for (PS) the following statements are equivalent:
	
	(i)  For any $(x_1(0),\cdots, x_N(0))$ satisfying A1),
	\begin{equation*}
	\sum_{i=1}^N\mathbb{E}\int_0^{\infty} e^{-\rho t} \left(\|\hat{x}_i(t)\|^2+\|\hat{u}_i(t)\|^2\right)dt<\infty.
	\end{equation*}	
	
	(ii) (\ref{eq15}) and (\ref{eq16}) admit $\rho$-stabilizing solutions\footnote{For a Riccati equation (\ref{eq15}), $P$ is called a $\rho$-stabilizing solution if $P$
		satisfies (\ref{eq15}) and all the eigenvalues of $A-BR^{-1}B^TP$ are in left half-plane.}, respectively, and $\bar{A}+G-\frac{\rho}{2}I$ is Hurwitz.
	
	(iii) A2) and A4) hold.
\end{theorem}

\begin{remark}
   $M_1$ and $M_2$ are Hamiltonian matrices. The Hamiltonian matrix plays a significant role in studying general algebraic Riccati equations. See more details of the property of Hamiltonian matrices in \cite{AFIJ03, M77}.
\end{remark}

To show Theorem \ref{thm7}, we need two lemmas. The first lemma is a result from \cite[Theorem 6]{M77}.
\begin{lemma}\label{lem5}
	Equations (\ref{eq15}) and (\ref{eq16}) admit $\rho$-stabilizing solutions if and only if
	A2) holds and both $M_1$ and $M_2$ have no eigenvalues on the imaginary axis.
\end{lemma}

\begin{lemma}\label{lem4}
	Let A1) hold. Assume that (\ref{eq15}) and (\ref{eq16}) admit $\rho$-stabilizing solutions, respectively, and $\bar{A}+G-\frac{\rho}{2}I$ is Hurwitz. Then
	\begin{equation*}
	\sum_{i=1}^N\mathbb{E}\int_0^{\infty} e^{-\rho t} \left(\|\hat{x}_i(t)\|^2+\|\hat{u}_i(t)\|^2\right)dt<\infty.
	\end{equation*}
\end{lemma}
\emph{Proof.} From the definition of $\rho$-stabilizing solutions, $A-BR^{-1}B^TP-\frac{\rho}{2}I$ and $A+G-BR^{-1}B^T\Pi-\frac{\rho}{2}I$ are Hurwitz. By the argument in the proof of Theorem \ref{thm4}, the lemma follows. $\hfill \Box$

\emph{The Proof of Theorem \ref{thm7}.}  By using Lemmas \ref{lem5} and \ref{lem4} together with a similar argument in the proof of Theorem \ref{thm4}, 
the Theorem follows. $\hfill \Box$

\begin{example}\label{ex1}
	Consider a scalar system with $A=a$, $B=b$, $G=g$, $Q=q$, $\Gamma=\gamma$, $R=r>0$. Then
	$$M_1=\left[\begin{array}{cc}
	a-{\rho}/{2} & b^2/r \\
	q& -a+{\rho}/{2}
	\end{array}\right],\quad
	M_2=\left[\begin{array}{cc}
	a+g-{\rho}/{2} & {b^2}/{r} \\
	q(1-\gamma)^2& -(a+g-{\rho}/{2})
	\end{array}\right].$$
	By direct computations, neither $M_1$ nor $M_2$ has eigenvalues in imaginary axis if and only if
	\begin{align}\label{eq41a}
	&(a-\frac{\rho}{2})^2+\frac{b^2}{r}q>0,\\ \label{eq42a}
	&(a+g-\frac{\rho}{2})^2+\frac{b^2}{r}(1-\gamma)^2q>0.
	\end{align}
	Note that if $q>0$ (or $a-{\rho}/{2}<0$,\ $q=0$), i.e., $(a-{\rho}/{2},\sqrt{q})$ is observable (detectable), then (\ref{eq41a}) holds,
	and  if $(1-\gamma)^2q>0$ ($a+g-{\rho}/{2}<0,$\ $q=0$), i.e., $(a+g-{\rho}/{2},\sqrt{q}(1-\gamma))$ is observable (detectable), then (\ref{eq42a}) holds.
	
	For this model, the Riccati equation (\ref{eq15}) is written as
	\begin{equation}\label{eq43a}
	\frac{b^2}{r}p^2-(2a-\rho)p-q=0.
	\end{equation}
	Let $\Delta=4[(a-{\rho}/{2})^2+{b^2q}/{r}]$. If (\ref{eq41a}) holds then $\Delta>0$, which implies (\ref{eq43a}) admits two solutions. If $q>0$ then (\ref{eq43a}) has a unique positive solution such that $a-b^2p/r-{\rho}/{2}=-\sqrt{\Delta}/2<0$.
	If $q=0$ and $a-\rho/2<0$ then (\ref{eq43a}) has a unique non-negative solution $p=0$ such that $a-b^2p/r-{\rho}/{2}=a-{\rho}/2<0$.
	
	Assume that (\ref{eq41a}) and (\ref{eq42a}) hold. By Theorem \ref{thm7}, the system is uniformly stable if and only if
	$(a-\rho/2,b)$ is stabilizable (i.e., $b\not=0$ or $a-\rho/2<0$), and $a-b^2p/r-{\rho}/{2}+g<0$. Note that $a-b^2p/r-{\rho}/{2}<0$. When $g\leq0$,
we have $a-b^2p/r-{\rho}/{2}+g<0$.
	
\end{example}

\begin{example}
	We further consider the model in Example \ref{ex1} for the case that $a+g=\rho/2$ and $\gamma=1$ (i.e., (\ref{eq42a}) does not hold).  In this case, the Riccati equation (\ref{eq16}) admits a unique solution $\Pi=0$.  (\ref{eq17}) becomes
	$\rho s=\dot{s}+\frac{\rho}{2}s$
	and has a unique solution $s=0$ in $C_{\rho/2}([0,\infty),\mathbb{R})$. Thus, $\bar{x}$ satisfies
	\begin{equation}\label{eq44}
	\frac{d\bar{x}}{dt}=\frac{\rho}{2}\bar{x}+f.
	\end{equation}
	Assume that $f$ is a constant.
	Then (\ref{eq44}) does not admit a solution in $C_{\rho/2}([0,\infty),\mathbb{R})$ unless $\bar{x}(0)=-{2f}/{\rho}$.
	
\end{example}

We are in a position to state the asymptotic optimality of
the decentralized control.

\begin{theorem}\label{thm8}
	Let A1)-A4) hold. For Problem (PS), the set of decentralized control laws
	$\{\hat{u}_1,\cdots,\hat{u}_N\}$ given by (\ref{eq14}) has asymptotic social optimality, i.e.,
	$$\Big|\frac{1}{N}J_{\rm soc}(\hat{u})-\frac{1}{N}\inf_{u\in \mathcal{U}_c}J_{\rm soc}(u)\Big|=O(\frac{1}{\sqrt{N}}).$$
	
\end{theorem}

\emph{Proof.} We first prove that for $u\in \mathcal{U}_c$, $J_{\rm soc}(u)< \infty$ implies that
\begin{equation}\label{eq36a}
\mathbb{E}\int_0^{\infty}e^{-\rho t}(\|x_i\|^2+\|u_i\|^2)dt<\infty,
\end{equation} for all $i=1,\cdots,N$.
From $J_{\rm soc}(u)< \infty$, we have
$\mathbb{E}\int_0^{\infty}e^{-\rho t}\|u_i\|^2dt<\infty$ and
\begin{equation}\label{eq36}
\mathbb{E}\int_0^{\infty}e^{-\rho t}\big\|x_i-\Gamma x^{(N)}\big\|^2_Qdt<\infty,
\end{equation}
which further implies that
\begin{equation}\label{eq37}
\begin{aligned}
&\mathbb{E}\int_0^{\infty}e^{-\rho t}\big\|(I-\Gamma)x^{(N)}\big\|^2_Q
\leq \frac{1}{N}\sum_{i=1}^N
\mathbb{E}\int_0^{\infty}e^{-\rho t}\big\|x_i-\Gamma x^{(N)}\big\|^2_Q dt<\infty.
\end{aligned}
\end{equation}
By (\ref{eq1}) we have
\begin{equation*}
\begin{aligned}
dx^{(N)}(t)=\ & \left[(A+G)x^{(N)}(t)+Bu^{(N)}(t)+f(t)\right]dt+\frac{1}{N}\sum_{i=1}^N\sigma(t) dW_i(t),
\end{aligned}
\end{equation*}
which leads to for any $r\in [0,1]$,
\begin{equation}\label{eq39}
\begin{aligned}
x^{(N)}(t)=\ &e^{(A+G)r}x^{(N)}(t-r)+\int_{t-r}^te^{(A+G)(t-\tau)}[Bu^{(N)}(\tau)+f(\tau)]d\tau\\
&+\frac{1}{N}\sum_{i=1}^N\int_{t-r}^te^{(A+G)(t-\tau)}\sigma(\tau) dW_i(\tau).
\end{aligned}
\end{equation}
By $J_{\rm soc}(u)< \infty$ and basic SDE estimates, we can find a constant $C$ such that
$$\mathbb{E}\int_r^{\infty}e^{-\rho t}\Big\|\int_{t-r}^te^{(A+G)(t-\tau)}Bu^{(N)}(\tau)d\tau\Big\|^2dt\leq C.$$
From (\ref{eq37}) and (\ref{eq39}) we obtain
$$
\begin{aligned}
\mathbb{E}\int_r^{\infty}&e^{-\rho t}[x^{(N)}(t-r)]^Te^{(A+G)^Tr}(I-\Gamma)^TQ(I-\Gamma)\cdot e^{(A+G)r}x^{(N)}(t-r)dt\leq C,
\end{aligned}
$$
which together with A3)
implies that
\begin{equation}\label{eq41}
\mathbb{E}\int_0^{\infty}e^{-\rho t}\|x^{(N)}(t)\|^2dt<\infty.
\end{equation}
This and (\ref{eq36}) lead to
\begin{equation}\label{eq41b}
\mathbb{E}\int_0^{\infty}e^{-\rho t}\|x_i(t)\|^2_Qdt<\infty.
\end{equation}

By (\ref{eq1}), we have
\begin{equation}\label{eq42}
\begin{aligned}
x_i(t)=\ &e^{Ar}x_i(t-r)+\int_{t-r}^te^{A(t-\tau)}[Bu_i(\tau)+f(\tau)+Gx^{(N)}(\tau)]d\tau\\
&+\int_{t-r}^te^{A(t-\tau)}\sigma(\tau) dW_i(\tau).
\end{aligned}
\end{equation}
It follows from (\ref{eq41}) that
\begin{align*}
&\mathbb{E}\int_r^{\infty}e^{-\rho t}\Big\|\int_{t-r}^te^{A(t-\tau)}Gx^{(N)}(\tau)d\tau\Big\|^2dt\cr
\leq & \mathbb{E}\int_0^{\infty}e^{-\rho \tau}\|Gx^{(N)}(\tau)\|^2\int_{0}^r\Big\|e^{(A-\frac{\rho}{2}I)v}\Big\|^2dvd\tau\leq C.
\end{align*}
From (\ref{eq41b}) and (\ref{eq42}), we obtain that
$$\mathbb{E}\int_r^{\infty}e^{-\rho t}x_i^T(t-r)e^{A^Tr}Qe^{Ar}x_i(t-r)dt\leq C.$$
This together with A3)
implies that
\begin{equation*}
\mathbb{E}\int_0^{\infty}e^{-\rho t}\|x_i(t)\|^2dt<\infty,
\end{equation*}
which gives (\ref{eq36a}). By Theorem \ref{thm4},
\begin{equation*}
\mathbb{E}\int_0^{\infty}e^{-\rho t}\big(\|\tilde{x}_i\|^2+\|\tilde{u}_i\|^2\big)dt<\infty.
\end{equation*}

By a similar argument to the proof of Theorem \ref{thm3} combined with Lemma \ref{lem2}, the conclusion follows.
$\hfill \Box$

If A3) is replaced by A3$^{\prime}$), the decentralized control (\ref{eq14}) still has asymptotic social optimality.
\begin{corollary}
	Assume that A1)-A2), A3$^{\prime}$), A4) hold.
	The set of decentralized control laws
	given by (\ref{eq14}) is asymptotically socially optimal.
\end{corollary}
\emph{Proof.}   Without loss of generality, we simply assume $A+G=\hbox{diag}\{\mathbb{A}_{1},\mathbb{A}_2\}$, where $\mathbb{A}_1-(\rho/2) I$ is Hurwitz, and $-(\mathbb{A}_2-(\rho/2) I)$ is Hurwitz (If necessary, we may apply a nonsingular linear transformation as in the proof of Theorem \ref{thm6}). Write $x^{(N)}=[z_1^T,z_2^T]$ and ${\hat{Q}}^{1/2}=[S_1,S_2]$ such that
$$\big\|(I-\Gamma)x^{(N)}\big\|_Q^2=
\|S_1{z}_1+S_2z_2\|^2,$$
and $(\mathbb{A}_2-(\rho/2) I,S_2)$ is observable which is due to the detectability of $(A+G-(\rho/2) I,\hat{Q}^{1/2})$. By the proof of Theorem \ref{thm4} or \cite{H10}, $\mathbb{E}\int_0^{\infty}e^{-\rho t}\|u^{(N)}\|^2dt<\infty$ implies $\mathbb{E}\int_0^{\infty}e^{-\rho t}\|z_1\|^2dt<\infty$, which together with (\ref{eq37}) gives  $\mathbb{E}\int_0^{\infty}e^{-\rho t}\|S_2z_2\|^2dt<\infty$. This and the observability of $(A_2-(\rho/2) I,S_2)$ leads to $\mathbb{E}\int_0^{\infty}e^{-\rho t}\|z_2\|^2dt<\infty$. Thus, $\mathbb{E}\int_0^{\infty}e^{-\rho t}\|x^{(N)}\|^2dt<\infty$.  The other parts of the proof are similar to that of Theorem \ref{thm8}.
$\hfill \Box$

\section{Mean Field LQ Games}

In this section, we investigate the game problem for LQ mean field systems.

\textbf{(PG)}. Seek a set of decentralized control laws to minimize individual cost
for each agent in the system (\ref{eq1})-(\ref{eq2}).
\subsection{The finite-horizon problem}\label{sec5.1}

We first consider the finite-horizon problem. 
Suppose that $\bar{x}\in C([0,T], \mathbb{R}^n)$ is given for approximation of $x^{(N)}$.
Replacing $x^{(N)}$ in (\ref{eq1}) and (\ref{eq3}) by $\bar{x}$, we have the following auxiliary optimal control problem.

$$\hbox{\textbf{(P2)}} \inf_{u_i\in L^2_{{\cal F}_t^i}(0, T; \mathbb{R}^r)} \bar{J}_i^{\rm F}(u_i),$$
where
\begin{equation*}\label{eq1g}
\begin{aligned}
d\grave{x}_i(t) = \ & [A\grave{x}_i(t)+Bu_i(t)+G\bar{x}(t)+f(t)]dt+\sigma(t) dW_i(t), \quad  1\leq i\leq N,\cr
\bar{J}_i^{\rm F}(u_i)=\ &\mathbb{E}\int_0^{T}e^{-\rho t}\Big \{\|\grave{x}_i-\Gamma\bar{x}-\eta\|^2_Q
	+\|u_i\|_R^2\Big\}dt.
\end{aligned}
\end{equation*}

By examining the variation of $\bar{J}_{i}^{\rm F}$, we obtain the unique optimal control of (P2).

\begin{theorem}
	Assume $Q\geq0, R>0$. Then the FBSDE
	\begin{equation}\label{eq47a}
	\left\{\begin{aligned}
	d\grave{x}_i &= \big[A\grave{x}_i-BR^{-1}B^Tp_i+G\bar{x}+f\big]dt+\sigma dW_i,\quad\\
	d{p}_i&=-\Big[
	(A^T-\rho I)p_i+
	Q{\grave{x}}_i-Q\Gamma\bar{x}-Q\eta\Big]dt+
	q_i dW_i, \\
	\grave{x}_i(&0)=x_{i0},\ p_i(T)=0,i=1,2,\cdots,N
	\end{aligned}\right.
	\end{equation}
	admits a unique solution $(\grave{x}_i,p_i,q_i)$, and the optimal control $\hat{u}_i=-R^{-1}B^Tp_i$.
\end{theorem}

\emph{Proof.} Since $Q\geq0$ and $R>0$, then by \cite{YZ99}, (P2) is uniformly convex, and hence admits a unique optimal control. By a similar argument with Theorem \ref{thm1}, the conclusion follows.
$\hfill \Box$


It follows  from (\ref{eq47a}) that

$$
\left\{\begin{aligned}
d\grave{x}^{(N)}= & \big[A\grave{x}^{(N)}-BR^{-1}B^Tp^{(N)}+G\bar{x}+f\big]dt+\frac{1}{N}\sum_{i=1}^N\sigma dW_i,\\
d{p}^{(N)}=&-\Big[
(A-\rho I)^Tp^{(N)}+
Q\grave{x}^{(N)}-Q\Gamma\bar{x}-Q\eta\Big]dt+\frac{1}{N}\sum_{i=1}^Nq_i dW_i,\\
\grave{x}^{(N)}(0)&=\frac{1}{N}\sum_{i=1}^Nx_{i0},\ p^{(N)}(T)=0.
\end{aligned}\right.$$
Replacing $\grave{x}^{(N)}$ by $\bar{x}$, we have
\begin{equation}\label{eq48}
\left\{\begin{aligned}
d\bar{x} &= \big[(A+G)\bar{x}-BR^{-1}B^T\bar{p}+f\big]dt,\quad  \bar{x}(0)=x_{0},\\
d\bar{p}&=-\big[
(A-\rho I)^T\bar{p}+
Q\bar{x}-Q\Gamma\bar{x}-Q\eta\big]dt, \  \bar{p}(T)=0.
\end{aligned}\right.
\end{equation}
Let $\bar{p}=\bar{P}\bar{x}+\hat{s}$. By It\^{o}'s formula, we obtain
\begin{align*}
&\bar{P}\left\{(A+G)\bar{x}-BR^{-1}B^T(\bar{P}\bar{x}+\hat{s})+f\right\}dt+(\dot{\bar{P}}\bar{x}+\dot{\hat{s}})dt
\cr
&=d\bar{p}
=-\big[(A-\rho I)^T(\bar{P}\bar{x}+\hat{s}) +Q\bar{x}-Q\Gamma\bar{x}-Q\eta\big]dt.
\end{align*}
This implies
\begin{align}\label{eq8b}
&\rho \bar{P}=\dot{\bar{P}}+{A}^T\bar{P}+\bar{P}(A+G)-\bar{P}BR^{-1}B^T\bar{P}+Q-Q\Gamma,\ \bar{P}(T)=0,\\
&\rho \hat{s}=\dot{\hat{s}}+(A-BR^{-1}B^T\bar{P})^T\hat{s}+\bar{P}f-Q\eta,\  \hat{s}(T)=0. \label{eq9b}
\end{align}

Denote $\tilde{p}_i=p_i-\bar{p}$, and $\tilde{x}_i=\grave{x}_i-\bar{x}$. Then by (\ref{eq47a}) and (\ref{eq48}) we have
$$
\left\{\begin{aligned}
d\tilde{x}_i &= \big[A\tilde{x}_i-BR^{-1}B^T\tilde{p}_i\big]dt+\sigma dW_i,  \tilde{x}_i(0)=x_{i0}-\bar{x}_{0},\\
d\tilde{p}_i&=-\big[
(A-\rho I)^T\tilde{p}_i+
Q\tilde{x}_i\big]dt+q_i dW_i, \  \tilde{p}_i(T)=0.
\end{aligned}\right.$$
Let $\tilde{p}_i=P\tilde{x}_i$. By It\^{o}'s formula,
$$
\begin{aligned}
d\tilde{p}_i=\ &-\big[
(A-\rho I)^T\tilde{p}_i+
Q\tilde{x}_i\big]dt+q_i dW_i\\
=\ &\dot{P}\tilde{x}_idt+P[(A\tilde{x}_i-BR^{-1}B^TP\tilde{x}_i)dt+\sigma dW_i],
\end{aligned}$$
which implies that $q_i=P\sigma$, and
\begin{equation}\label{eq10c}
\rho P=\dot{P}+A^TP+PA-PBR^{-1}B^TP+Q.
\end{equation}
Assume

\textbf{A5)} Equation (\ref{eq8b}) admits a solution in $C([0,T],\mathbb{R}^n)$.

By the local Lipschitz-continuous property of the quadratic function, (\ref{eq8b}) can admit a unique local solution in a small time duration $[T_0, T]$.
It may be referred to \cite{AFIJ03} for some sufficient conditions of the existence of the solution in $[0,T]$.
We now provide a necessary and sufficient condition to guarantee the global solvability of (\ref{eq8b}).

\begin{proposition}
	(\ref{eq8b})
	admits a solution in $C([0,T],\mathbb{R}^n)$
if and only if
	for any $t \in [0, T],$
	$$
	\det\{(0, I)e^{\mathcal{A}t}(0, I)^{T}\}>0,
	$$
	where
	$$\mathcal{A}=\left(\begin{array}{cc} A+G& -BR^{-1}_2B^T \\
	Q\Gamma-Q & -(A-\rho I)^{T}
	\end{array}\right).$$
\end{proposition}
\emph{Proof.} Sufficiency is given by \cite[Theorem 4.3, p.48]{MY99}. Necessity is implied from Proposition 4.2 and Theorem 3.2 of \cite[Chapter 2]{MY99}.
$\hfill \Box$

Let
\begin{equation}\label{eq12c}
\hat{u}_i=-R^{-1}B^T[P\hat{x}_i+(\bar{P}-P)\bar{x}+\hat{s}],
\end{equation}
where $P, \bar{P}$ and $\hat{s}$ are determined by (\ref{eq10c}), (\ref{eq8b}) and (\ref{eq9b}), respectively, and $\bar{x}$ and $\hat{x}_i$ satisfy
\begin{align}
\label{eq53}
d\bar{x} =\  &\big[A\bar{x}-BR^{-1}B^T(\bar{P}\bar{x}+\hat{s})+G\bar{x}+f\big]dt,\  \bar{x}(0)=x_{0},\\
\label{eq53a}
d\hat{x}_i =\ & \big[(A-BR^{-1}B^TP)\hat{x}_i-BR^{-1}B^T[(\bar{P}-P)\bar{x}+\hat{s}]\cr
&+G\hat{x}^{(N)}+f\big]dt+\sigma dW_i,\ \hat{x}_i(0)=x_{i0}.
\end{align}

Denote $u_{-i}=(u_1,\cdots,u_{i-1},u_{i+1},\cdots,u_N)$.
\begin{theorem}\label{thm10}
	Let A1), A5) hold and $Q\geq 0$. The set of decentralized strategies
	$\{\hat{u}_1,\cdots,\hat{u}_N\}$ given by (\ref{eq12c}) is an $\varepsilon$-Nash equilibrium, i.e.,
	\begin{equation}\label{eq54}
	\inf_{u_i\in L^2_{{\cal F}_t}(0, T; \mathbb{R}^r)}J_i^{\rm F}(u_i,\hat{u}_{-i})\geq J_i^{\rm F}(\hat{u}_i,\hat{u}_{-i})-\varepsilon,
	\end{equation}
	where $\varepsilon=({1}/{\sqrt{N}})$. 
\end{theorem}

\emph{Proof.} See the Appendix \ref{app d}. $\hfill \Box$

\subsection{The infinite-horizon problem}
For simplicity, we consider the case $G=0$.

Based on the analysis in Section \ref{sec5.1}, we may design the following decentralized control for (PG):
\begin{equation}\label{eq14c}
\begin{aligned}
\hat{u}_i(t)= &-{R^{-1}}B^T[P\hat{x}_i(t)+(\bar{P}-P)\bar{x}(t)+\hat{s}(t)],\ t\geq 0, \ i=1,\cdots, N,
\end{aligned}
\end{equation}
where $P$ and $\bar{P}$ are determined by
\begin{align}\label{eq15c}
\rho P=&A^TP+PA-PBR^{-1}B^TP+Q,\\
\label{eq16c}
\rho \bar{P}= &A^T\bar{P}+\bar{P} A-\bar{P} BR^{-1}B^T\bar{P}+Q(I-\Gamma),
\end{align}
respectively, and
$\hat{s}, \bar{x}\in C_{\rho/2}([0,\infty),\mathbb{R}^n)$ are determined by
\begin{align}
\rho \hat{s}=&\dot{\hat{s}}+[A-BR^{-1}B^T\bar{P}]^T\hat{s}+\bar{P} f-{\eta},\label{eq17c}\\
\frac{d\bar{x}}{dt}=&A\bar{x}-B{R^{-1}}B^T(\bar{P}\bar{x}+\hat{s})+f,\   \bar{x}(0)=\bar{x}_0.\label{eq18c}
\end{align}
and $\hat{x}_i$ satisfies
\begin{eqnarray}
  d\hat{x}_i \hspace*{-0.1cm}&\hspace*{-0.1cm}=\hspace*{-0.1cm}&\hspace*{-0.1cm} \big[(A-BR^{-1}B^TP)\hat{x}_i-BR^{-1}B^T[(\bar{P}-P)\bar{x}+\hat{s}]\cr
&\hspace*{-0.2cm}&+G\hat{x}^{(N)}+f\big]dt+\sigma dW_i,\ \hat{x}_i(0)=x_{i0}.\label{eq17d}
\end{eqnarray}
Here the existence conditions of $P, \bar{P},s$ and $\bar{x}$ need to be investigated further.

We introduce the following assumptions.

\textbf{A6)} $(A-\frac{\rho}{2}I, B)$ is stabilizable, $Q\geq0$ and $(A-\frac{\rho}{2}I, \sqrt{Q})$ is detectable.

\textbf{A7)} (\ref{eq16c}) admits a stabilizing solution. 

\begin{lemma}\label{lem2c}
	Assume that $M_3$ has $n$ stable eigenvalues (with negative real parts) and $n$ unstable eigenvalues,
	where $$M_3=\left[\begin{array}{cc}
	A-\frac{\rho}{2}I & BR^{-1}B^T \\
	Q(I-\Gamma)& -A^{T}+\frac{\rho}{2}I
	\end{array}\right].$$ Suppose that
	\begin{equation}\label{eq70}
	M_3\left[\begin{array}{c}
	L_1 \\
	L_2
	\end{array}\right]=\left[\begin{array}{c}
	L_1 \\
	L_2
	\end{array}\right]H_{11},
	\end{equation}
	where $H_{11}$ is Hurwitz and $L_1$ is invertible.
	Then A7) holds.
\end{lemma}

\emph{Proof.}
 Let $\bar{P}=-L_2L_1^{-1}$. It follows from (\ref{eq70}) that
\begin{equation}\label{eq71}
M_3\left[\begin{array}{c}
-I \\
\bar{P}
\end{array}\right]=\left[\begin{array}{c}
-I \\
\bar{P}
\end{array}\right] L_1H_{11} L_1^{-1}.
\end{equation}
By pre-multiplying by $[\bar{P} \ \ I]$ on both sides, we obtain
$$[\bar{P} \ \ I]
M_3\left[\begin{array}{c}
-I \\
\bar{P}
\end{array}\right]=0,
$$
which leads to (\ref{eq16c}). By (\ref{eq71}), we have
$A-BR^{-1}B^T\bar{P}-\frac{\rho}{2}I=L_1H_{11} L_1^{-1}$ is Hurwitz.
It is straightforward that
$s, \bar{x}\in C_{\rho/2}([0,\infty),\mathbb{R}^n)$. $\hfill \Box$

\begin{remark}
	The above lemma provides a convenient method to compute the stabilizing solutions of algebraic Riccati equations.
	Assume there exists an invertible matrix
	$V=\left[\begin{array}{cc}
	V_{11} & V_{12} \\
	V_{21}& V_{22}
	\end{array}\right]$ such that
	$V^{-1}M_3 V=\left[\begin{array}{cc}
	H_{11} & H_{12} \\
	0& H_{22}
	\end{array}\right],$
	where $V_{11}$ is invertible, and $H_{11}, -H_{22}$ are Hurwitz.
	Then $V_{21}V_{11}^{-1}$ is the stabilizing solution of (\ref{eq16c}). $V$ comprises $2n$ independent vectors, which are called Schur vectors \cite{L79}.
\end{remark}
\begin{lemma}\label{lem7}
	Assume that A1), A6), A7) hold.
	Then (\ref{eq17c})-(\ref{eq18c}) admit a set of unique solutions $s, \bar{x}\in C_{\rho/2}([0,\infty),\mathbb{R}^n)$,
and
\begin{equation*}
	\sum_{i=1}^N\mathbb{E}\int_0^{\infty} e^{-\rho t} \left(\|\hat{x}_i(t)\|^2+\|\hat{u}_i(t)\|^2\right)dt<\infty.
	\end{equation*}	
\end{lemma}

\emph{Proof.} By a similar argument in the proof of Theorem \ref{thm6}, 
the lemma follows. $\hfill \Box$

%
%
%

\begin{theorem}\label{thm11}
	Let A1), A6), A7) hold. For Problem (PG), the set of decentralized strategies
	$\{\hat{u}_1,\cdots,\hat{u}_N\}$ given by (\ref{eq14c}) is an $\varepsilon$-Nash equilibrium, i.e.,
	\begin{equation*}
	\inf_{u_i\in \mathcal{U}_{c,i}}J_i(u_i,\hat{u}_{-i})\geq J_i(\hat{u}_i,\hat{u}_{-i})-\varepsilon,
	\end{equation*}
	where $\varepsilon=({1}/{\sqrt{N}}).$
\end{theorem}

\emph{Proof.} See Appendix \ref{app d}.  $\hfill \Box$

\section{Comparison of Different Solutions}

In this section, we compare the proposed decentralized control laws 
with the feedback decentralized strategies in previous works \cite{HCM07, HCM12}.

We first introduce a definition from \cite{BO82}.
\begin{definition}
	For a control problem with an admissible control set $\mathcal{U}$, a control law $u\in \mathcal{U}$ is said to be a representation of another control
	$u^*\in \mathcal{U}$ if
	
	(i) they both generate the same unique state trajectory, and
	
	(ii) they both have the same open-loop value on this trajectory.
\end{definition}

For Problem (PS), let $f=0$, and $G=0$.
In \cite[Theorem 4.3]{HCM12}, the decentralized control laws are given by
\begin{equation}\label{eq47}
\breve{u}_i=-R^{-1}B^T(Px_i+\bar{s}),\quad i=1,\cdots,N,
\end{equation}
where $P$ is the semi-positive definite solution of (\ref{eq15c}), and $\bar{s}=\bar{K}{x}^{\dag}+\phi.$ Here $\bar{K}$ satisfies
\begin{align*}
  \rho\bar{K}=\ &\bar{K}\bar{A}+\bar{A}^T\bar{K}
-\bar{K}BR^{-1}B^T\bar{K}^T-Q_{\Gamma},
\end{align*}
and ${x}^{\dag}, \phi\in C_{\rho/2}([0,\infty),\mathbb{R}^n) $ are determined by
\begin{align*}
\frac{d\bar{x}^\dag}{dt}=\ &\bar{A}\bar{x}^{\dag}-BR^{-1}B^T(\bar{K}\bar{x}^\dag+\phi),\bar{x}^\dag(0)=\bar{x}_0,\cr
\frac{d\phi}{dt}=\ &-[A-BR^{-1}B^T(P+\bar{K})-{\rho}I]\phi+\bar{\eta},
\end{align*}
in which $\bar{A}=A-BR^{-1}B^TP$.
By comparing this with (\ref{eq16})-(\ref{eq18}), one can obtain that $\bar{K}=\Pi-P$, $\bar{x}=\bar{x}^{\dag}$ and $\phi=s$. From the above discussion, we have the equivalence of the two sets of decentralized control laws.
\begin{proposition}
	The set of decentralized control laws $\{\hat{u}_1,\cdots,\hat{u}_N\}$ in (\ref{eq14}) is
	a representation of $\{\breve{u}_1,\cdots,\breve{u}_N\}$ given by (\ref{eq47}).
\end{proposition}

For Problem (PG), let $f=0$, and $G=0$.
In \cite{HCM07}, the decentralized strategies are given by
\begin{equation}\label{eq74}
{u}_i^*=-R^{-1}B^T(Px_i+{s}^*),\quad i=1,\cdots,N,
\end{equation}
where $P$ is the positive definite solution of (\ref{eq15}), ${s}^*$ is determined by the fixed-point equation
\begin{equation}\label{eq75}
\left\{  \begin{aligned}
\rho s^*&=\frac{ds^*}{dt}+\bar{A}^Ts^*-Q\Gamma(\bar{x}^*+\eta),\\ 
\frac{d\bar{x}^*}{dt}&=\bar{A}\bar{x}^*-BR^{-1}B^Ts^*, \ \bar{x}^*(0)= \bar{x}_0.
\end{aligned}
\right.
\end{equation}

We now show the equivalence of the decentralized open-loop and feedback solutions to mean field games.
\begin{proposition}
	The set of decentralized control laws $\{\hat{u}_1,\cdots,\hat{u}_N\}$ in (\ref{eq14c}) is
	a representation of $\{{u}_1^*,\cdots,{u}_N^*\}$ given by (\ref{eq74}).
\end{proposition}

\emph{Proof.} Let $s^*=K^*\bar{x}^*+\psi$. 
From (\ref{eq75}), we have
\begin{align*}
\frac{ds^*}{dt}=&K^*\frac{d\bar{x}^*}{dt}+\frac{d\psi}{dt}\cr
=&K^*\big[\bar{A}\bar{x}^*-BR^{-1}B^T(K^*\bar{x}^*+\psi)\big]+\frac{d\psi}{dt}\cr
=&(\rho I-\bar{A})^T(K^*\bar{x}^*+\psi)+Q(\Gamma\bar{x}^*+\eta),
\end{align*}
which gives
\begin{align*}
\rho K^*&=K^*\bar{A}+\bar{A}^TK^*-K^*BR^{-1}B^TK^*-Q\Gamma,\\
\rho \psi&=\frac{d\psi}{dt}+[\bar{A}-BR^{-1}B^TK^*]^T \psi-Q\eta.
\end{align*}
By comparing this with (\ref{eq15c})-(\ref{eq17c}),
one can obtain $K=\bar{P}-P$, and $\psi=\hat{s}$. Thus, we have $u_i^*\equiv \hat{u}_i$, $i=1,\cdots,N,$
which implies that $\{{u}_1^*,\cdots,{u}_N^*\}$ is
a representation of  $\{\hat{u}_1,\cdots,\hat{u}_N\}$ in (\ref{eq14c}). $\hfill \Box$


\section{Numerical Examples}
In this section, some numerical examples are given to illustrate the effectiveness of the proposed decentralized control laws.

\begin{figure}[H]
	\centering
	\includegraphics[width=0.5\linewidth]{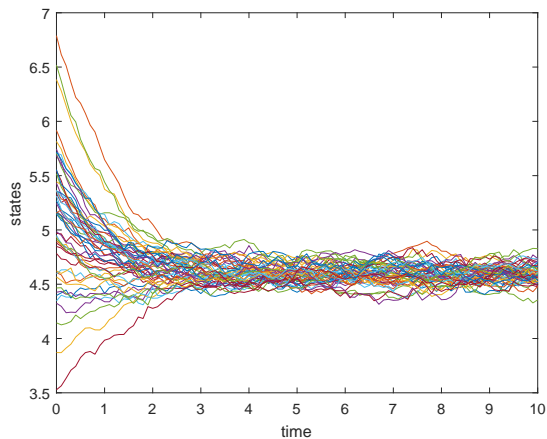}
	\caption{Curves of 50 agents with $A=0.2$.}
	\label{psa02g-02}
\end{figure}

\begin{figure}[H]
	\centering
	\includegraphics[width=0.5\linewidth]{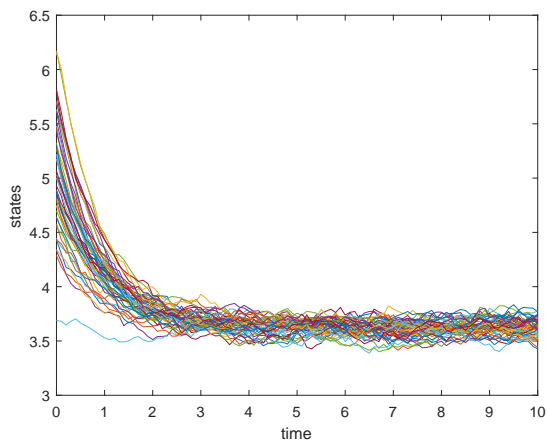}
	\caption{Curves of 50 agents with $A=1$.}
	\label{psa1g-02}
\end{figure}

We first consider a scalar system with $50$ agents in Problem (PS). Take $B=Q=R= 1,G=-0.2,f(t)=1,\sigma(t)=0.1, \rho=0.6, \Gamma=-0.2,\eta=5 $ in (\ref{eq1})-(\ref{eq2}). The initial states of $50$ agents are taken independently from a normal distribution $N(5,0.5)$.  Then, under the control law (\ref{eq14}), the state trajectories of agents for the cases with $ A=0.2$ and $A=1$ are shown in Figs. \ref{psa02g-02} and \ref{psa1g-02}, respectively. After the transient phase, the states of agents behave similarly and achieve agreement roughly.

\begin{figure}[H]
	\centering
	\includegraphics[width=0.55\linewidth]{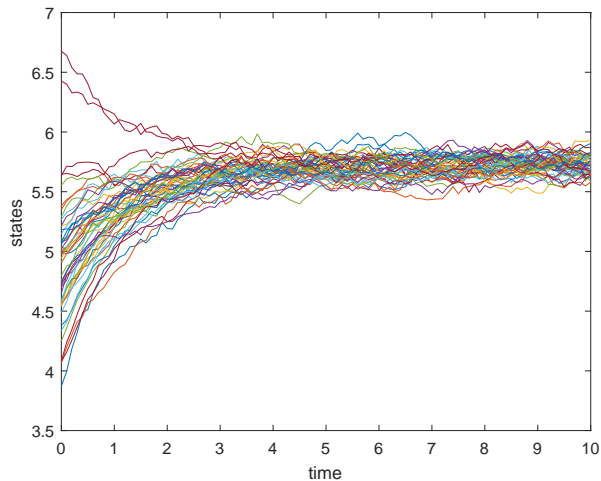}
	\caption{Curves of 50 agents with $A=0.2$.}
	\label{pga02g0}
\end{figure}

\begin{figure}[H]
	\centering
	\includegraphics[width=0.55\linewidth]{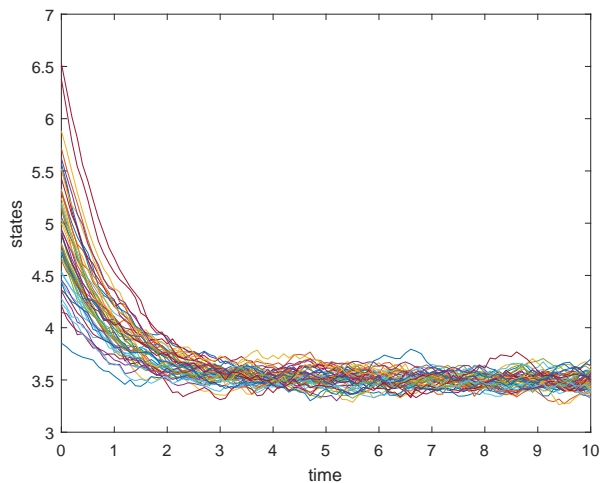}
	\caption{Curves of 50 agents with $A=1$.}
	\label{pga1g0}
\end{figure}

\begin{figure}[H]
	\centering
	\includegraphics[width=0.55\linewidth]{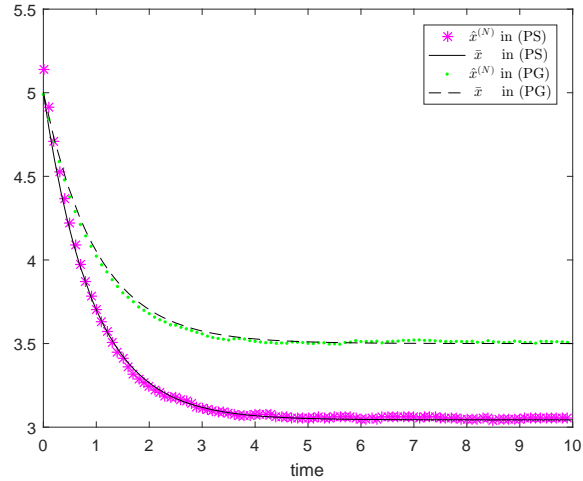}
	\caption{Curves of $\bar{x}$ and $\hat{x}^{(N)}$ in (PS) and (PG).}
	\label{ps-pg-a1g0}
\end{figure}

Next, we simulate the scalar case of Problem (PG), where the parameters are the same as above, except $G=0$. After the control laws (\ref{eq14c}) are applied, the state trajectories of 50 agents with $ A=0.2$ and $A=1$ are shown in Figs. \ref{pga02g0} and \ref{pga1g0}, respectively.
For the case $A=1$ and $G=0$, the trajectories of $\bar{x}$ and $\hat{x}^{(N)}$ in Problems (PS) and (PG) are shown in Fig. \ref{ps-pg-a1g0}.
It can be seen that $\bar{x}$ and $\hat{x}^{(N)}$ coincide well, which illustrate the consistency of mean field approximations.
Clearly, the state average of agents has significantly lower value in Problem (PS) than in (PG).

\begin{figure}[H]
	\centering
	\includegraphics[width=0.55\linewidth]{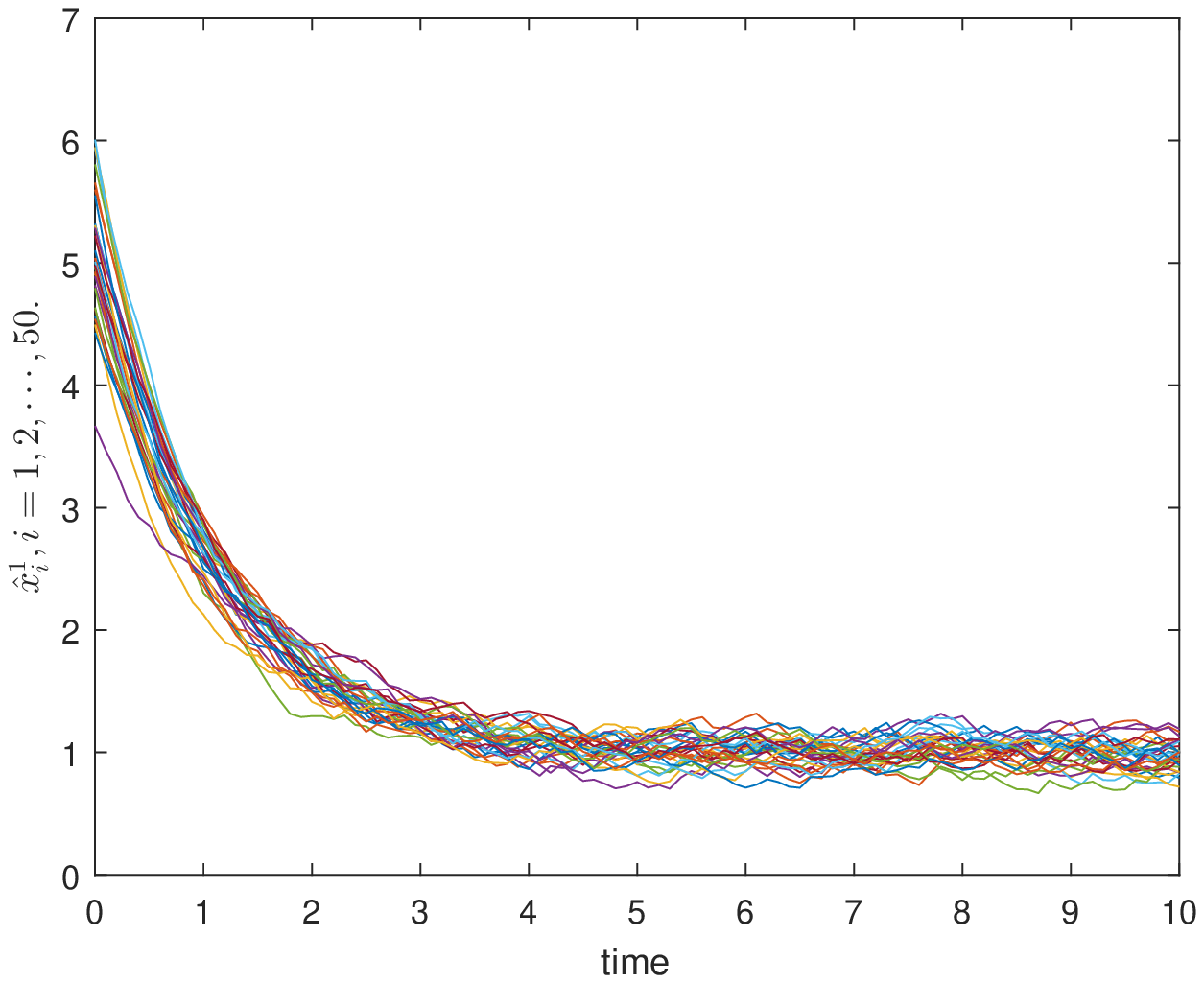}
	\caption{Curves of $\hat{x}^1_i$, $i=1,\cdots,50$.}
	\label{vector1}
\end{figure}

\begin{figure}[H]
	\centering
	\includegraphics[width=0.55\linewidth]{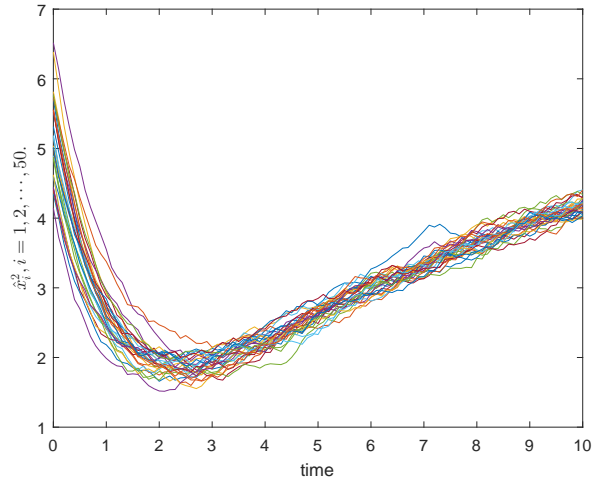}
	\caption{Curves of $\hat{x}^2_i$, $i=1,\cdots,50$.}
	\label{vector2}
\end{figure}

Finally, we consider the 2-dimensional case of Problem (PS). Take parameters as follows:
$A = \left[\begin{array}{cc}
0.1 & 0\\
-1  & 0.2\\
\end{array}\right]$,
$B = \left[\begin{array}{cc}
1 & 0\\
0  & 1\\
\end{array}\right]$,
$G = \left[\begin{array}{cc}
-0.5 & 0\\
0  & -0.3\\
\end{array}\right]$,
$B = \left[\begin{array}{c}
1 \\
1\\
\end{array}\right]$,
$Q = \left[\begin{array}{cc}
1 & 0\\
0  & 1\\
\end{array}\right]$,
$\Gamma = \left[\begin{array}{cc}
1 & 0\\
1  & 1\\
\end{array}\right]$,
$R = \left[\begin{array}{cc}
1 & 0\\
0  & 1\\
\end{array}\right]$,
$\eta = \left[\begin{array}{c}
0\\
0.5\\
\end{array}\right]$, $f = [1\ \ 1]^T$and $\sigma = [0.5\ \ 0.5]^T$. Denote
$\hat{x}_i(t)=[\hat{x}^1_i(t)\ \hat{x}^2_i(t)]^T$.
 Both of $\hat{x}^1_i(0)$ and $\hat{x}^2_i(0)$ are taken independently from a normal distribution $N(5,0.5)$. Under the control laws (\ref{eq14}), the trajectories of $\hat{x}^1_i$ and $\hat{x}^2_i$, $i=1,\cdots,N$ are shown in Figs. \ref{vector1} and \ref{vector2}, respectively.


\section{Concluding Remarks}

In this paper, we have considered uniform stabilization and asymptotic optimality for mean field LQ multiagent systems. 
For social control and Nash game problems, we design the decentralized open-loop control laws by the variational analysis, respectively, which are further shown to be asymptotically optimal. Two equivalent conditions are further given for uniform stabilization of the systems in different cases. Finally, we show such decentralized control laws are equivalent to the feedback strategies in previous works.

An interesting generalization is to consider mean field LQ control systems with partial measurements by using variational analysis. Also, the variational analysis may be applied to general nonlinear model to construct decentralized control laws for social control and Nash games.

\appendices
\section{Proof of Theorem \ref{thm3}}\label{app a}
\def\theequation{A.\arabic{equation}}
\setcounter{equation}{0}

To prove Theorem \ref{thm3}, we need a lemma.

\begin{lemma}\label{lem1}
	Let A1) hold and $Q\geq 0$. Under the control (\ref{eq15a}), we have
	\begin{equation}
	\max_{0\leq t\leq T}\mathbb{E}\|\hat{x}^{(N)}(t)-\bar{x}(t)\|^2=O(\frac{1}{N}).
	\end{equation}	
\end{lemma}
\emph{Proof.} It follows by (\ref{eq20}) that
\begin{equation*}
\begin{aligned}
d\hat{x}^{(N)}=\ &\big[(\bar{A}+G)\hat{x}^{(N)}-BR^{-1}B^T(K\bar{x}+s)+f\big]dt+\frac{1}{N}\sum_{i=1}^N\sigma dW_i.
\end{aligned}
\end{equation*}
From this and (\ref{eq12a}), we have
\begin{equation*}
d(\hat{x}^{(N)}-\bar{x})=(\bar{A}+G)(\hat{x}^{(N)}-\bar{x})dt+\frac{1}{N}\sum_{i=1}^N\sigma dW_i,
\end{equation*}
which leads to
\begin{equation}\label{eq21}
\begin{aligned}
\hat{x}^{(N)}(t)-\bar{x}(t)=\ &e^{(\bar{A}+G)t}[\hat{x}^{(N)}(0)-\bar{x}(0)]+\frac{1}{N}\sum_{i=1}^N\int_0^te^{(\bar{A}+G)(t-\tau)}\sigma dW_i(\tau).
\end{aligned}
\end{equation}
By A1), one can obtain
\begin{equation}
\begin{aligned}
&\mathbb{E}\big\|  \hat{x}^{(N)}(t)-\bar{x}(t)\big\|^2\\
\leq\ &\big\|2e^{(\bar{A}+G)t}\big\|^2\Big\{\mathbb{E}\big\|\hat{x}^{(N)}(0)-\bar{x}(0)\big\|^2+\frac{1}{N}\int_0^ttr\big[\sigma^T
e^{-(\bar{A}^T+G^T+\bar{A}+G)\tau}\sigma\big] d\tau\Big\}\cr
\leq \ &  \frac{2}{N}\big\|e^{(\bar{A}+G)t}\big\|^2\Big\{\max_{1\leq i\leq N}\mathbb{E}\|\hat{x}_{i0}\|^2+\int_0^Ttr\big[\sigma^T
e^{-(\bar{A}^T+G^T+\bar{A}+G)\tau}\sigma\big] d\tau\Big\},
\end{aligned}
\end{equation}
which completes the proof.  $\hfill \Box$

\emph{Proof of Theorem \ref{thm3}}.
We first prove that for $u\in \mathcal{U}_c$,  $J_{\rm soc}^{\rm F}(u)< \infty$ implies that
$\mathbb{E}\int_0^Te^{-\rho t}(\|x_i\|^2+\|u_i\|^2)dt<\infty,$ for all $i=1,\cdots,N$.
By $J_{\rm soc}^{\rm F}(u)< \infty$, we have $\mathbb{E}\int_0^Te^{-\rho t}\|u_i\|^2dt<\infty.$ This leads to
$$\mathbb{E}\int_0^Te^{-\rho t}\|u^{(N)}\|^2dt\leq\frac{1}{N}\sum_{i=1}^N\mathbb{E}\int_0^Te^{-\rho t}\|u_i\|^2dt<\infty,$$
where $u^{(N)}=\frac{1}{N}\sum_{i=1}^Nu_i.$
By (\ref{eq1}),
\begin{equation*}\label{eq38}
\begin{aligned}
dx^{(N)}(t)=\ & \left[(A+G)x^{(N)}(t)+Bu^{(N)}(t)+f(t)\right]dt+\frac{1}{N}\sum_{i=1}^N\sigma(t) dW_i(t),
\end{aligned}
\end{equation*}
which with A1) implies that
$$\max_{0\leq t\leq T}\mathbb{E}\|{x}^{(N)}(t)\|^2\leq C.$$
Note that
$$x_i(t)=e^{At}x_{i0}+\int_0^te^{A(t-\tau)}[Gx^{(N)}(\tau)+Bu_i(\tau)+f(\tau)]d\tau.$$
We have
\begin{equation}\label{eq22a}
\begin{aligned}
&\mathbb{E}\int_0^Te^{-\rho t}\|x_i\|^2dt
\leq C\Big(\mathbb{E}\|x_{i0}\|^2+\max_{0\leq t\leq T}\mathbb{E}\|{x}^{(N)}(t)\|^2+\max_{0\leq t\leq T}\mathbb{E}\|u_i(t)\|^2\Big)<\infty.
\end{aligned}
\end{equation}
By (\ref{eq12a}) and (\ref{eq20}), we obtain that
\begin{equation}\label{eq22b}
\mathbb{E}\int_0^Te^{-\rho t}\big(\|\hat{x}_i\|^2+\|\hat{u}_i\|^2+\|\bar{x}\|^2)dt<\infty.
\end{equation}

Let $\tilde{x}_i=x_i-\hat{x}_i$,   $\tilde{u}_i=u_i-\hat{u}_i$ and $\tilde{x}^{(N)}=\frac{1}{N}\sum_{i=1}^N \tilde{x}_i$. Then by (\ref{eq1}) and (\ref{eq20}),
\begin{equation}\label{eq32}
d\tilde{x}_i=(A\tilde{x}_i+{G}\tilde{x}^{(N)}+B\tilde{u}_i)dt, \  \tilde{x}_i(0)=0.
\end{equation}
From (\ref{eq3}), we have
\begin{equation}
\begin{aligned}
J^{\rm F}_{\rm soc}(u)=\ &\sum_{i=1}^N\mathbb{E}\int_0^Te^{-\rho t}\Big[\big\|\hat{x}_i-\Gamma \hat{x}^{(N)}-\eta +\tilde{x}_i-\Gamma \tilde{x}^{(N)}\big\|^2_Q+\big\|\hat{u}_i+\tilde{u}_i\big\|^2_{R}\Big]dt
\\
=\ &\sum_{i=1}^N(J_i^{\rm F}(\hat{u})+\tilde{J}_i^{\rm F}(\tilde{u})+I_i),
\end{aligned}
\end{equation}
where
\begin{align*}
\tilde{J}_i^{\rm F}(\tilde{u})\stackrel{\Delta}{=}\mathbb{E}\int_0^Te^{-\rho t}\big[&\|\tilde{x}_i-\Gamma \tilde{x}^{(N)}\|^2_Q+\|\tilde{u}_i\|^2_{R}\big]dt,\\
I_i=2\mathbb{E}\int_0^Te^{-\rho t}\Big[&\big(\hat{x}_i-\Gamma \hat{x}^{(N)}-\eta\big)^TQ\big(\tilde{x}_i-\Gamma \tilde{x}^{(N)}\big)+\hat{u}_i^TR\tilde{u}_i\Big]dt.
\end{align*}
By (A4), $\tilde{J}_i^{\rm F}(\tilde{u})\geq 0$. We now prove $\frac{1}{N}\sum_{i=1}^N I_i=O(\frac{1}{\sqrt{N}})$.
\begin{equation}\label{eq24a}
\begin{aligned}
&\sum_{i=1}^N I_i\\
=\ &\sum_{i=1}^N 2\mathbb{E}\int_{0}^{T}e^{-\rho t}\Big\{\tilde{x}_i^T\big[Q(\hat x_i-\Gamma \hat{x}^{(N)}-\eta)-\Gamma^TQ((I-\Gamma)\hat{x}^{(N)}-\eta)\big]+\sum_{i=1}^N\hat{u}_i^TR\tilde{u}_i\Big\}dt\cr
=\ & \sum_{i=1}^N 2\mathbb{E}\int_{0}^{T}e^{-\rho t}\Big\{\tilde{x}_i^T\big[Q(\hat x_i-\Gamma\bar{x}-\eta)-\Gamma^TQ((I-\Gamma)\bar{x}-\eta)\big]+\sum_{i=1}^N\hat{u}_i^TR\tilde{u}_i\Big\}dt\cr
&+\sum_{i=1}^N2\mathbb{E}\int_{0}^{T}e^{-\rho t}(\hat{x}^{(N)}-\bar{x})^TQ_{\Gamma}\tilde{x}_idt.
\end{aligned}
\end{equation}
By (\ref{eq8a})-(\ref{eq10a}), (\ref{eq32}) and It\^{o}'s formula,
$$  \begin{aligned}
0=\ &\sum_{i=1}^N\mathbb{E}\big[ e^{-\rho T} \tilde{x}_i^T(T)(P\hat{x}_i(T)+K\bar{x}(T)+s(T))-  \tilde{x}_i^T(0)(P\hat{x}_i(0)+K\bar{x}(0)+s(0))\big]\cr
=\ &\mathbb{E}\int_{0}^{T}\sum_{i=1}^Ne^{-\rho t}\Big\{- \tilde{x}_i^T\big[ Q\hat{x}_i-Q(\Gamma \bar{x}+\eta) 
-\Gamma^TQ\left((I-\Gamma) \bar{x}-\eta\right)\big]
\-\hat{u}_i^TR\tilde{u}_i)\Big\}dt\\&+N\mathbb{E}\int_{0}^{T}e^{-\rho t}(\hat{x}^{(N)}-\bar{x})^T(G^TP+PG)\tilde{x}^{(N)}dt.
\end{aligned}$$
From this and (\ref{eq24a}), we obtain
$$
\begin{aligned}
\frac{1}{N}\sum_{i=1}^N I_i=2\mathbb{E}\int_{0}^{T}&e^{-\rho t}(\hat{x}^{(N)}-\bar{x})^T\cdot (Q_{\Gamma}+G^TP+PG)\tilde{x}^{(N)}dt.
\end{aligned}
$$
By Lemma \ref{lem1}, (\ref{eq22a}) and (\ref{eq22b}), we obtain
$$
\begin{aligned}
\Big|\frac{1}{N}\sum_{i=1}^N I_i\Big|^2\leq C &\mathbb{E}\int_{0}^{T}e^{-\rho t}\|\hat{x}^{(N)}-\bar{x}\|^2dt\cdot \mathbb{E}\int_{0}^{T}e^{-\rho t}\|\tilde{x}^{(N)}\|^2dt,
\end{aligned}
$$
which implies $|\frac{1}{N}\sum_{i=1}^N I_i|=O(1/\sqrt{N})$.  $\hfill \Box$

\section{Proofs of Lemma \ref{lem2} and Theorem \ref{thm4}}\label{app b}
\def\theequation{B.\arabic{equation}}
\setcounter{equation}{0}

\emph{Proof of Lemma \ref{lem2}.} 
From (\ref{eq21}), we have
\begin{equation*}\label{eq21b}
\begin{aligned}
\hat{x}^{(N)}(t)-\bar{x}(t)=\ &e^{(\bar{A}+G)t}[\hat{x}^{(N)}(0)-\bar{x}(0)]+\frac{1}{N}\sum_{i=1}^N\int_0^te^{(\bar{A}+G)(t-v)}\sigma dW_i(v).
\end{aligned}
\end{equation*}
Thus,
$$\begin{aligned}
&\mathbb{E}\int_0^{\infty} e^{-\rho t} \left(\|\hat{x}^{(N)}(t)-\bar{x}(t)\|^2\right)dt \\
\leq\ &  2\mathbb{E}\int_0^{\infty} \left\| e^{(\bar{A}+G-\frac{\rho}{2}I)t}\right\|^2\big\|\hat{x}^{(N)}(0)-\bar{x}(0)\big\|^2dt\cr
&+2\mathbb{E}\int_0^{\infty}e^{-\rho t} \frac{1}{N}  \left\| \int_0^te^{(\bar{A}+G)(t-v)}\sigma dW_i(v)\right\|^2dt\cr
\leq\ & 2\int_0^{\infty} \left\| e^{(\bar{A}+G-\frac{\rho}{2}I)t}\right\|^2\mathbb{E}\big\|\hat{x}^{(N)}(0)-\bar{x}(0)\big\|^2dt\cr
&+\frac{2}{N} \mathbb{E}\int_0^{\infty}e^{-\rho t}\int_0^t tr\left[\sigma^T \sigma e^{(\bar{A}+G+\bar{A}^T+G^T)(t-v)}\right] dvdt\cr
\leq\ & \frac{2}{N}\int_0^{\infty} \left\| e^{(\bar{A}+G-\frac{\rho}{2}I)t}\right\|^2\mathbb{E}\big\|\max_{1\leq i\leq N}\hat{x}_i(0)\big\|^2dt\cr
&+ \frac{C}{N} \mathbb{E}\int_0^{\infty}e^{-\rho v}\|\sigma\|^2\int_v^{\infty}\big\|e^{(\bar{A}+\bar{G}-\frac{\rho}{2}I)(t-v)}\big\|^2 dtdv\cr
\leq\ & O(\frac{1}{N}).
\end{aligned}$$
$\hfill \Box$

\emph{Proof of Theorem \ref{thm4}.} By A1)-A4), Lemmas \ref{lem2a} and \ref{lem2}, we obtain that $\bar{x}\in C_{\rho/2}([0,\infty),\mathbb{R}^n)$ and
$$\mathbb{E}\int_0^{\infty} e^{-\rho t} \left(\big\|\hat{x}^{(N)}(t)-\bar{x}(t)\big\|^2\right)dt=O(\frac{1}{N}),$$
which further gives that 
$$\mathbb{E}\int_0^{\infty} e^{-\rho t} \|\hat{x}^{(N)}(t)\|^2dt<\infty.$$
Denote $g\stackrel{\Delta}{=}-BR^{-1}B^T((\Pi-P)\bar{x}+s)+Gx^{(N)}+f$. Then $\mathbb{E}\int_0^{\infty} e^{-\rho t} \|g(t)\|^2dt<\infty$  and 
\begin{equation}\label{eq14a}
\hat{x}_i(t)=e^{\bar{A}t}\hat{x}_{i0}+\int_0^te^{\bar{A}(t-v)}g(v)dv+\int_0^te^{\bar{A}(t-v)}\sigma dW_i.
\end{equation}
Note that $\bar{A}-\frac{\rho}{2}I$ is Hurwitz. By Schwarz's inequality,
$$\begin{aligned}
&\mathbb{E}\int_0^{\infty}e^{-\rho t}\|\hat{x}_i(t)\|^2dt\\
\leq \ &3 \mathbb{E}\int_0^{\infty} \left\| e^{(\bar{A}-\frac{\rho}{2}I)t}\right\|^2\|\hat{x}_{i0}\|^2dt+3\mathbb{E}\int_0^{\infty}e^{-\rho t}t\int_0^t\Big\|e^{\bar{A}(t-v)}g(v)\Big\|^2dvdt\cr
&+3 \mathbb{E}\int_0^{\infty}e^{-\rho t}\int_0^t tr[e^{\bar{A}^T(t-v)}\sigma^T(v)\sigma(v)e^{\bar{A}(t-v)}] dvdt\cr
\leq \ &C+3\mathbb{E}\int_0^{\infty}e^{-\rho v}\|g(v)\|^2\int_v^{\infty}t\big \|e^{(\bar{A}-\frac{\rho}{2}I)(t-v)}\big\|^2dtdv\cr
&+3 C\mathbb{E}\int_0^{\infty}e^{-\rho v}\|\sigma(v)\|^2\int_v^{\infty}\big \|e^{(\bar{A}-\frac{\rho}{2}I)(t-v)}\big\|^2dtdv\cr
\leq \ &C+3C\mathbb{E}\int_0^{\infty}e^{-\rho v}\|g(v)\|^2dv+3 C\mathbb{E}\int_0^{\infty}e^{-\rho v}\|\sigma(v)\|^2dv\leq C_1.
\end{aligned}$$
This with (\ref{eq14}) completes the proof.   \hfill $\Box$

\section{Proofs of Theorems \ref{thm5} and \ref{thm6}}\label{app c}
\def\theequation{C.\arabic{equation}}
\setcounter{equation}{0}

\emph{Proof.}
i)$\Rightarrow$ ii). 
By (\ref{eq20}),
\begin{equation}\label{eq24}
\begin{aligned}
\frac{d\mathbb{E}[\hat{x}_i]}{dt}=\ &\bar{A}\mathbb{E}[\hat{x}_i]-BR^{-1}B^T((\Pi-P)\bar{x}+s)+G\mathbb{E}[\hat{x}^{(N)}]+f, \quad \mathbb{E}[\hat{x}_i(0)]=\bar{x}_0.
\end{aligned}
\end{equation}
It follows from A1) that
$$\mathbb{E}[\hat{x}_i]=\mathbb{E}[\hat{x}_j]=\mathbb{E}[\hat{x}^{(N)}], \ j\not =i.$$
By comparing (\ref{eq18}) and (\ref{eq24}), we obtain $\mathbb{E}[\hat{x}_i]=\bar{x}$. Note that $\|\bar{x}\|^2=\big\|\mathbb{E}\hat{x}_i\big\|^2\leq \mathbb{E}\|\hat{x}_i\|^2$.
It follows from (\ref{eq23}) that
\begin{equation}\label{eq25}
\int_0^{\infty}e^{-\rho t} \|\bar{x}(t)\|^2dt<\infty.
\end{equation}
By (\ref{eq18}), we have
$$
\begin{aligned}
\bar{x}(t)=&e^{(A+G-BR^{-1}B^T\Pi)t}\Big[\bar{x}_0+\int_0^te^{-(A+G-BR^{-1}B^T\Pi)\tau}h(\tau)d\tau\Big] ,
\end{aligned}
$$
where $h=-BR^{-1}B^Ts+f$. By the arbitrariness of $\bar{x}_0$ with (\ref{eq25}) we obtain that $A+G-BR^{-1}B^T\Pi-\frac{\rho}{2}I$ is Hurwitz. That is, $(A+G-\frac{\rho}{2}I, B)$ is stabilizable.
By \cite{AM90}, (\ref{eq16}) admits a unique solution such that $\Pi>0$. Note that $\mathbb{E}[x^{(N)}]^2\leq \frac{1}{N}\sum_{i=1}^N\mathbb{E}[\hat{x}_i^2]$. Then from (\ref{eq23}) we have 
\begin{equation}\label{eq27}
\mathbb{E} \int_0^{\infty}e^{-\rho t}\big\|\hat{x}^{(N)}(t)\big\|^2dt<\infty.
\end{equation}
This leads to $ \mathbb{E} \int_0^{\infty}e^{-\rho t}\|g(t)\|^2dt<\infty$, where $g{=}-BR^{-1}B^T((\Pi-P)\bar{x}+s)+G\hat{x}^{(N)}+f$.
By (\ref{eq14a}), we obtain 
$$
\begin{aligned}
\mathbb{E}\|\hat{x}_i(t)\|^2=\ & \mathbb{E}\left\|e^{\bar{A}t}\left(x_{i0}+\int_0^te^{-\bar{A}\tau}g(\tau)d\tau\right)\right\|^2 +\mathbb{E}\int_0^ttr\big[\sigma^T(\tau)e^{(\bar{A}^T+\bar{A})(t-\tau)}
\sigma(\tau)\big]d\tau.
\end{aligned}
$$
By (\ref{eq23}) and the arbitrariness of ${x}_{i0}$  we obtain that $\bar{A}-\frac{\rho}{2}I$ is Hurwitz, i.e., $(A-\frac{\rho}{2}I, B)$ is stabilizable.
By \cite{AM90}, (\ref{eq15}) admits a unique solution such that $P>0$.

From (\ref{eq25}) and (\ref{eq27}),
\begin{equation}\label{eq28}
\mathbb{E} \int_0^{\infty}e^{-\rho t} \big\|\hat{x}^{(N)}(t)-\bar{x}(t)\big\|^2 dt<\infty.
\end{equation}
On the other hand, (\ref{eq21}) gives
\begin{equation*}
\begin{aligned}
&\mathbb{E}\big\|\hat{x}^{(N)}(t)-\bar{x}(t)\big\|^2\\
=\ &\mathbb{E}\big\|e^{(\bar{A}+G)t}[\hat{x}^{(N)}(0)-\bar{x}_0]\big\|^2+\frac{1}{N}\int_0^ttr\big[\sigma^T(\tau)e^{(\bar{A}^T+G^T+\bar{A}+G)(t-\tau)}\sigma(\tau)\big] d\tau.
\end{aligned}
\end{equation*}
By (\ref{eq28}) and the arbitrariness of ${x}_{i0}, i=1,\cdots,N$, we obtain that $\bar{A}+G-\frac{\rho}{2}I$ is Hurwitz.

(ii)$\Rightarrow$(iii).  Define $V(t)=e^{-\rho t}\bar{y}^T(t)\Pi \bar{y}(t)$,
where $\bar{y}$ satisfies
\begin{equation*}
\frac{d\bar{y}}{dt}=(A+G)\bar{y}+B\bar{u},\quad  \bar{y}(0)=\bar{y}_0.
\end{equation*}
Denote $V$ by $V^*$ when $\bar{u}=\bar{u}^*=-{R^{-1}}B^T\Pi \bar{y}$. By (\ref{eq16}) we have
\begin{align*}
\frac{dV^*}{dt}=\ &\bar{y}^T(t)\big[-\rho\Pi+(A+G-B{R^{-1}}B^T\Pi)^T\Pi+\Pi(A+G-B{R^{-1}}B^T\Pi) \big]\bar{y}(t)\cr
=\ &\bar{y}^T(t)\big[-\hat{Q}-\Pi B{R^{-1}}B^T\Pi\big]\bar{y}(t)\leq 0.
\end{align*}
Note that $V^*\geq0$. Then $\lim_{t\to\infty}V^*(t)$ exists, which implies
\begin{equation}\label{eq43}
\lim_{t_0\to\infty}[V^*(t_0)-V^*(t_0+T)]=0.
\end{equation}

Rewrite $\Pi(t)$ in (\ref{eq11}) by $\Pi_{T}(t)$. Then we have $\Pi_{T+t_0}(t_0)=\Pi_{T}(0)$.
By (\ref{eq11}),
\begin{align*}
&\int_{t_0}^{T+t_0}e^{-\rho t} (\bar{y}^{T}\hat{Q}\bar{y}+\bar{u}^TR\bar{u})dt\cr
=\ &e^{-\rho t_0}\bar{y}^T({t_0})\Pi_{T+t_0}(t_0)\bar{y}({t_0})+ \int_0^Te^{-\rho t} \big\|\bar{u}+{R^{-1}}B^T\Pi_{T+t_0}(t_0) \bar{y}\big\|^2_Rdt\cr
\geq\ &e^{-\rho t_0}\big\|\bar{y}({t_0})\big\|^2_{\Pi_{T+t_0}(t_0)} =e^{-\rho t_0}\big\|\bar{y}({t_0})\big\|^2_{\Pi_{T}(0)}.
\end{align*}
This with (\ref{eq43}) implies
\begin{align*}
&\lim_{t_0\to\infty}e^{-\rho t_0}\big\|\bar{y}({t_0})\big\|^2_{\Pi_{T}(0)}\cr
\leq&\lim_{t_0\to\infty} \int_{t_0}^{T+t_0}e^{-\rho t} (\|\bar{y}\|_{\hat{Q}}^2+\|\bar{u}^*\|^2_R)dt
=\lim_{t_0\to\infty}[V^*(t_0)-V^*(t_0+T)]=0.
\end{align*}
By A3), one can obtain that there exists $T>0$ such that $\Pi_{T}(0)>0$ (See e.g. \cite{ZQ16, ZZC08}).
Thus, we have $\lim_{t\to\infty}e^{-\rho t}\big\|\bar{y}({t})\big\|^2=0$, which $(A+G-\frac{\rho}{2} I, B)$ is stabilizable. Similarly, we can show $(A-\frac{\rho}{2} I, B)$ is stabilizable.

(iii)$\Rightarrow$(i).
This part has been proved in Theorem \ref{thm4}.  $\hfill \Box$

\emph{Proof of Theorem \ref{thm6}.}
(iii)$\Rightarrow$(i). From \cite{AM90}, (\ref{eq15}) and (\ref{eq16}) admit unique solutions $P\geq0, \Pi\geq0$ such that $A-BR^{-1}B^TP-\frac{\rho}{2}I$ and $A-BR^{-1}B^T\Pi-\frac{\rho}{2}I$
are Hurwitz, respectively. Thus, there exists a unique $s(0)$ such that $s\in C_{\rho/2}([0,\infty),\mathbb{R}^n)$. It is straightforward that $\bar{x}\in C_{\rho/2}([0,\infty),\mathbb{R}^n)$. By the argument in the proof of Theorem \ref{thm4}, (i) follows.

(i)$\Rightarrow$(ii). The proof of this part is similar to that of (i)$\Rightarrow$(ii) in Theorem \ref{thm5}.

(ii)$\Rightarrow$(iii). Since $\Pi\geq0$, then there exists an orthogonal $U$ such that
$$U^T\Pi U=\left[
\begin{array}{cc}
0 & 0 \\
0& \Pi_{2}
\end{array}
\right],$$
where $\Pi_2>0$.
From (\ref{eq15}),
\begin{align}\label{eq40}
\rho U^T\Pi U=&(U^T\bar{\mathbb{A}}U)^TU^T\Pi U+U^T\Pi UU^T\bar{\mathbb{A}}U+U^T\bar{Q}U,
\end{align}
where $\bar{\mathbb{A}}\stackrel{\Delta}{=}A+G-\Pi BR^{-1}B^T\Pi,\bar{Q}=\hat{Q}+\Pi BR^{-1}B^T\Pi $.
Denote
$$U^T\bar{\mathbb{A}}U=\left[\begin{array}{cc}
\bar{\mathbb{A}}_{11}& \bar{\mathbb{A}}_{12} \\
\bar{\mathbb{A}}_{21}& \bar{\mathbb{A}}_{22}
\end{array}\right], \ U^T\bar{Q}U=\left[\begin{array}{cc}
\bar{Q}_{11}& \bar{Q}_{12} \\
\bar{Q}_{21}& \bar{Q}_{22}
\end{array}\right].$$
By pre- and post-multiplying by $\xi^T$ and $\xi$ where $\xi=[\xi_1^T,0]^T$, it follows that
$$0=\rho\xi^T U^T\Pi U\xi=\xi^T U^T\bar{Q}U\xi.$$
From the arbitrariness of $\xi_1$, we obtain $\bar{Q}_{11}=0$.
Since $\bar{Q}$ is semi-positive definite, then $\bar{Q}_{12}=\bar{Q}_{21}=0$, and $\bar{Q}_{22}\geq0$. By comparing each block matrix of both sides of (\ref{eq40}),
we obtain $\bar{\mathbb{A}}_{21}=0$. It follows from (\ref{eq40}) that
\begin{equation}\label{eq40b}
\rho\Pi_2=\Pi_2\bar{\mathbb{A}}_{22}+\bar{\mathbb{A}}_{22}^T\Pi_2+\bar{Q}_{22}.
\end{equation}

Let $\zeta=[\zeta_1^T,\zeta_2^T]^T=U^T\bar{y}^*$, where $\bar{y}^*$ satisfies $\dot{\bar{y}}^*=\bar{\mathbb{A}}\bar{y}^*$.
Then we have
\begin{align*}
\dot{\zeta_1}&=\bar{\mathbb{A}}_{11}\zeta_1+\bar{\mathbb{A}}_{12}\zeta_2,\cr
\dot{\zeta_2}&=\bar{\mathbb{A}}_{22}\zeta_2.
\end{align*}
By Lemma 4.1 of \cite{W68}, the detectability of $(A+G, \hat{Q}^{1/2})$ implies the detectability of $(\bar{\mathbb{A}}, \bar{Q}^{1/2})$.
Take $\zeta(0)=\xi=[\xi_1^T,0]^T$. Then $\bar{Q}^{1/2}\bar{y}=\bar{Q}^{1/2}U\zeta=0$, which together with the detectability of $(\bar{\mathbb{A}}, \bar{Q}^{1/2})$ implies $\zeta_1\to 0$ and $\bar{\mathbb{A}}_{11}$ is Hurwitz.
Denote $S(t)=e^{-\rho t}\zeta_2^T\Pi_2\zeta_2$. By (\ref{eq40b}),
$$S(T)-S(0)=-\int_0^T\zeta_2(t)^T\bar{Q}_{22}\zeta_2(t)dt\leq 0,$$
which implies $\lim_{t\to\infty}S(t)$ exists.
By a similar argument with the proof of Theorem \ref{thm5}, we obtain $\lim_{t_0\to\infty}e^{-\rho t_0}\big\|\zeta_2(t_0)\big\|^2_{\Pi_{2,T}(0)}=0$
and $\Pi_{2,T}(0)>0$, which gives $\zeta_2\to 0$ and $\bar{\mathbb{A}}_{22}$ is Hurwitz. This with the fact that $\bar{\mathbb{A}}_{11}$ is Hurwitz gives that
$\zeta$ is stable, which leads to (iii). $\hfill \Box$

\section{Proof of Theorems \ref{thm10} and \ref{thm11}}\label{app d}
\def\theequation{D.\arabic{equation}}
\setcounter{equation}{0}

\emph{Proof of Theorem \ref{thm10}.}  
From (\ref{eq53}) and (\ref{eq53a}), we have
$$d(\hat{x}^{(N)}-\bar{x})=(\bar{A}+G)(\hat{x}^{(N)}-\bar{x})dt+\frac{\sigma}{N}\sum_{i=1}^NdW_i,$$
where $\bar{A}=A-BR^{-1}B^TP$. This implies that
\begin{equation}\label{eq55}
\sup_{0\leq t\leq T}\mathbb{E}\|\hat{x}^{(N)}(t)-\bar{x}(t)\|^2=O(\frac{1}{N}).
\end{equation}
By Schwarz's inequality,
\begin{align}\label{eq57}
J_i^{\rm F}(\hat{u}_i,\hat{u}_{-i})\leq \ & \bar{J}_i^{\rm F}(\hat{u}_i)+\mathbb{E}\int_0^Te^{-\rho t}\|\hat{x}^{(N)}(t)-\bar{x}(t)\|^2dt\cr
&+2C\Big(\mathbb{E}\int_0^Te^{-\rho t}\|\hat{x}^{(N)}(t)-\bar{x}(t)\|^2dt\Big)^{1/2}\cr
\leq\ &  \bar{J}_i^{\rm F}(\hat{u}_i)+O(1/\sqrt{N}).
\end{align}

To prove (\ref{eq54}), it suffices to only consider $u_i\in L^2_{{\cal F}_t}(0, T; \mathbb{R}^r)$ 
 such that $J_i^{\rm F}(u_i,\hat{u}_{-i})\leq J_i^{\rm F}(\hat{u}_i,\hat{u}_{-i})<\infty$.
By (\ref{eq3}),
\begin{equation}\label{eq58}
\mathbb{E}\int_0^Te^{-\rho t}\|u_i\|^2dt<\infty.
\end{equation}
After the set of strategies $(u_i,\hat{u}_{-i})$ is applied, the corresponding dynamics of $N$ agents can be written as
$$\begin{aligned}
dx_i=&(Ax_i+Bu_i+Gx^{(N)}+f)dt+\sigma dW_i,\\
dx_j=&(Ax_j+B\hat{u}_j+Gx^{(N)}+f)dt+\sigma dW_j,\ j=1,\cdots,i-1,i+1,\cdots,N.
\end{aligned}
$$
This with (\ref{eq48}) implies
\begin{align*}
&d(x^{(N)}-\bar{x})\\
=\ &[(A+G)(x^{(N)}-\bar{x})+\frac{1}{N}B(u_i-\hat{u}_i)+B(\hat{u}^{(N)}+R^{-1}B^T\bar{p})]dt+\frac{1}{N}\sum_{j=1}^N\sigma dW_j\cr
=\ &[(A+G)(x^{(N)}-\bar{x})+\frac{1}{N}B(u_i-\hat{u}_i)-BR^{-1}B^TP(\hat{x}^{(N)}-\bar{x})]dt+\frac{\sigma}{N}\sum_{j=1}^NdW_j.
\end{align*}
By (\ref{eq55}), (\ref{eq58}) and elementary SDE estimates, one can obtain
\begin{equation}\label{eq59}
\mathbb{E}\int_0^Te^{-\rho t}\|{x}^{(N)}-\bar{x}\|^2dt<O(\frac{1}{N}).
\end{equation}
We have
$$d(x_i-\grave{x}_i)=[A(x_i-\grave{x}_i)+G(x^{(N)}-\bar{x})]dt,$$
which together with (\ref{eq59}) gives that
\begin{equation}\label{eq60}
\mathbb{E}\int_0^Te^{-\rho t}\|x_i-\grave{x}_i\|^2dt<O(\frac{1}{N}).
\end{equation}
Note that
\begin{align*}
&\|x_i-(\Gamma x^{(N)}+\eta)\|^2_Q\geq \|\grave{x}_i-(\Gamma \bar{x}+\eta)\|^2_Q+2[\grave{x}_i-(\Gamma \bar{x}+\eta)]^TQ[(x_i-\grave{x}_i)+\Gamma(\bar{x}-x^{(N)})],
\end{align*}
and $\bar{J}_i^{\rm F}(u_i)<\infty$.
By Schwarz's inequality, (\ref{eq59}) and (\ref{eq60}), we obtain
\begin{align*}
&J_i^{\rm F}(u_i,\hat{u}_{-i})\cr
\geq\ &\bar{J}_i^{\rm F}(u_i)-\Big[\mathbb{E}\int_0^Te^{-\rho t} \|\grave{x}_i
-(\Gamma \bar{x}+\eta)\|^2_Qdt\Big]^{{1}/{2}}\cr
&\cdot \Big[\mathbb{E} \int_0^Te^{-\rho t} \|(x_i-\grave{x}_i)+\Gamma(\bar{x}-x^{(N)}\|^2_Qdt\Big]^{{1}/{2}}\cr
\geq\ &\bar{J}_i^{\rm F}(u_i)-O(1/\sqrt{N}).
\end{align*}
From this and (\ref{eq57}), the theorem follows.  $\hfill \Box$

\emph{Proof of Theorem \ref{thm11}.} Note that  $\{\hat{x}_i(t), i=1,\cdots,N\}$ are mutually independent processes with the expectation $\bar{x}(t)$. By Lemma \ref{lem7}, 
\begin{equation*}
\begin{aligned}
&\mathbb{E}\int_0^{\infty} e^{-\rho t} \|\hat{x}_i(t)-\bar{x}(t)\|^2dt
\leq \frac{1}{N}\mathbb{E}\int_0^{\infty} e^{-\rho t} \|\hat{x}_i(t)\|^2dt=O(\frac{1}{N}).
\end{aligned}
\end{equation*}
We only need to show $\mathbb{E}\int_0^{\infty}e^{-\rho t}\|{x}_i
\|^2_{Q}dt\leq C$ for all
$u_i$ satisfying
\begin{equation}\label{eq66}
J_i({u}_i, \hat{u}_{-i})\leq J_i(\hat{u}_i, \hat{u}_{-i})\leq C_0.
\end{equation}
From (\ref{eq66}), we obtain
$$\begin{aligned}
&\mathbb{E}\int_0^{\infty}e^{-\rho t}\big\|{x}_i
- \frac{1}{N}{x}_i\big\|^2_{Q}dt\cr
\leq\ & \mathbb{E}\int_0^{\infty}e^{-\rho t}\big\|{x}_i
- \frac{1}{N}{x}_i-\frac{1}{N}\sum_{j\not=i}\hat{x}_j\big\|^2_{Q}dt +\mathbb{E}\int_0^{\infty}e^{-\rho t}\big\|\frac{1}{N}\sum_{j\not=i}\hat{x}_j\big\|^2_{Q}dt\cr
\leq\  &C_0+ \mathbb{E}\int_0^{\infty}e^{-\rho t}\frac{1}{N-1}\sum_{j\not=i}\big\|\hat{x}_j\big\|^2_{Q}dt\leq C,
\end{aligned}$$
which with Lemma \ref{lem7} implies
\begin{equation*}
\mathbb{E}\int_0^{\infty}e^{-\rho t}\|{x}_i
\|^2_{Q}dt\leq C_1,
\end{equation*}
where $C_1$ is independent of $N$.
The rest of the proof follows by that of Theorem \ref{thm10}. $\hfill \Box$

\begin{IEEEbiography}{Bingchang Wang}
 received the M.Sc. degree in Mathematics from
Central South University, Changsha, China, in 2008, and
the Ph.D. degree in System Theory from Academy of
Mathematics and Systems Science, Chinese Academy of
Sciences, Beijing, China, in 2011. From September 2011
to August 2012, he was with Department of Electrical
and Computer Engineering, University of Alberta, Canada,
as a Postdoctoral Fellow. From September 2012 to September 2013, he was with School of Electrical Engineering and Computer
Science, University of Newcastle, Australia, as a Research
Academic.

From October 2013, he has
been with School of Control Science and Engineering, Shandong University, China, as an associate Professor. He held visiting
appointments as a Research Associate with Carleton University, Canada, from November 2014 to May 2015, and with the Hong Kong Polytechnic University from November 2016 to January
 2017. He also visited
 the Hong Kong Polytechnic University as a Research Fellow in March 2017 and May 2018.
His current research interests include mean field games, stochastic control, multiagent
systems and event based control. He received the IEEE CSS Beijing Chapter Young Author Prize in 2018.
\end{IEEEbiography}

\begin{IEEEbiography}{Huanshui Zhang}
(SM'06) received the B.S. degree
in mathematics from Qufu Normal University, Shandong,
China, in 1986, the M.Sc. degree in control
theory from Heilongjiang University, Harbin, China,
in 1991, and the Ph.D. degree in control theory from
Northeastern University, China, in 1997.

He was a Postdoctoral Fellow at Nanyang Technological
University, Singapore, from 1998 to 2001
and Research Fellow at Hong Kong Polytechnic
University, Hong Kong, China, from 2001 to 2003.
He is currently holds a Professorship at Shandong
University, Shandong, China. He was a Professor with the Harbin Institute
of Technology, Harbin, China, from 2003 to 2006. He also held visiting
appointments as a Research Scientist and Fellow with Nanyang Technological
University, Curtin University of Technology, and Hong Kong City University
from 2003 to 2006. His interests include optimal estimation and control,
time-delay systems, stochastic systems, signal processing and wireless sensor
networked systems.
\end{IEEEbiography}

\end{document}